\documentclass[11pt]{article}
\usepackage{amsmath,amssymb,epsf}
\usepackage[french]{babel}
\usepackage[applemac]{inputenc}
\advance \textwidth 3cm
\advance \textheight 3cm
\usepackage{amsfonts}
\usepackage{latexsym}
\newtheorem{theorem}{Th\'eor\`eme}
\newtheorem{lemme}{Lemme}
\newtheorem{corollaire}{Corollaire}

\newtheorem{prop}{Proprit}
\newtheorem{remarque}{Remarque}
\newenvironment{preuve}[1]{\par\noindent\underline{Preuve #1} :\quad}%
{\unskip\nobreak\hfil\penalty50\hskip2em\null\nobreak\hfil%
$\Box$\parfillskip0pt\par\medskip}

\newcommand{\trace}{\mathrm{Tr}}

\title{Valeur propre minimale d'une matrice Toeplitz et d'un produit de matrices de Toeplitz.}

\author{ Philippe Rambour\thanks{Universit\'{e} de Paris Sud,
      B\^atiment 425; F-91405
Orsay Cedex;
tel : 01 69 15 57 28 ; fax 01 69 15 60 19
      \mbox{e-mail : philippe.rambour@math.u-psud.fr}}
      }
\date{}
\begin{document}
\maketitle
  \renewcommand{\abstractname}{R\'ESUM\'E}
  \begin{abstract}
    \textbf{Valeur propre minimale d'une matrice de Toeplitz et d'un produit de matrices de Toeplitz.}\\
     Nous donnons une expression asymptotique 
     de la plus petite valeur propre $\lambda_{N,\alpha}$ de la matrice 
     $T_N(\varphi_{\alpha})$ o\`u 
          $\varphi_{\alpha}(e^{i \theta})=\vert 1- e^{i \theta} \vert ^{2\alpha} c_{1}(e^{i \theta})$, 
    avec $c_{1}$ une fonction strictement positive suffisamment r\'eguli\`ere
     et  $0<\alpha< \frac{1}{2}$. Nous obtenons $\lambda_{N,\alpha}\sim c_{\alpha}N^{-2\alpha}c_{1}(1)$ 
     et nous donnons un encadrement de $c_{\alpha}$.
     Pour obtenir un \'equivalent de la valeur propre minimale nous donnons et utilisons un thorme qui relie les coefficients de $T_N^{-1}(\varphi_\alpha)$ et 
     ceux de $T_N (\varphi^{-1}_\alpha)$. Sous l'hypothse $\alpha_1+\alpha_2 >\frac{1}{2}$ nous obtenons galement une expression asymptotique de la valeur propre minimale de $T_N (\varphi _{\alpha_1}) T_N (\varphi _{\alpha_2}) $.
     
            \end{abstract}
     \renewcommand{\abstractname}{ABSTRACT}
          \begin{abstract}
          \textbf{Minimal eigenvalue of a Toeplitz matrix and of a product of Toeplitz matrices.}\\
This paper is essentially devoted to the study of the minimal eigenvalue $\lambda_{N,\alpha}$  of the Toepllitz 
matrice $T_N(\varphi_{\alpha})$ where $\varphi_{\alpha}(e^{i \theta})=\vert 1- e^{i \theta} \vert ^{2\alpha} c_{1}(e^{i \theta})$ with
$c_{1}$ a positive sufficiently smooth function and $0<\alpha<\frac{1}{2}$. We obtain  
$\lambda_{N,\alpha}\sim c_{\alpha}N^{-2\alpha}c_{1}(1)$ when $N$ goes to the infinity and 
we have the bounds of $c_{\alpha}$. To obtain the asymptotic of $\lambda_{N,\alpha}$ we 
give a theorem which suggests that the entries of $T_N^{-1}(\varphi_{\alpha})$ and 
 $T_N (\varphi^{-1}_\alpha)$ are closely related. 
 If $\alpha_1 + \alpha_2 > \frac{1}{2}$ we obtain 
 the asymptotic of the minimal eigenvalue of 
 $T_N (\varphi _{\alpha_1}) T_N (\varphi _{\alpha_2}).$
 \end{abstract}


\textbf{\large{Mots clef}}

\textbf{Matrices de Toeplitz, produit de matrices de Toeplitz, valeur propre minimale, op\'erateurs \`a noyau.}
\section {Introduction}
 Rappelons  que si $f$ est une fonction de $L^1(\mathbb T)$
 on appelle matrice de Toeplitz d'ordre $N$ de symbole $f$, 
 et on note $T_N(f)$, la matrice $(N+1) \times (N+1)$ telle que 
 $\left(T_N (f)\right) _{k+1,l+1} =\hat f (l-k) \quad \forall \,k,l \quad 
 0\le k,l \le N$ o  $\hat h (j)$ dsigne le coefficient de Fourier 
 d'ordre $j$ d'une fonction $h$ (une bonne r\'ef\'erence peut \^{e}tre \cite{Bo.3}). Une fonction de  $L^1(\mathbb T)$ strictement positive sur le tore est appel\'ee une fonction r\'eguli\`ere. Dans ce travail on s'intresse  l'expression asymptotique de la valeur propre minimale des matrices de Toeplitz de symbole $\varphi_{\alpha}$ avec 
$\varphi_{\alpha}(e^{i \theta})=\vert 1- e^{i \theta} \vert ^{2\alpha} c_{1}(e^{i \theta})$, $c_{1}$
tant une fonction r\'eguli\`ere et aussi  la valeur propre minimale du produit de deux matrices de ce type. 
Ce travail complte et prolonge les articles \cite{BoVi}, \cite{BoW}, \cite{RS1111}et \cite{Ramb10}. 
 \\
Dans \cite{BoVi} B\"{o}ttcher et Virtanen donnent un \'equivalent quand $N$ tend vers l'infini de la valeur propre maximale de $T_{N} (\varphi_{\alpha})$ avec $-\frac{1}{2}<\alpha<0$. La norme $\Vert T_{N} \varphi_{\alpha}\Vert$ (qui est aussi la plus grande valeur propre) 
est reli\'ee \`a la norme d'un op\'erateur. On a (dans \cite{BoVi} on s'int\'eresse \`a des  $\alpha$ n\'egatifs)
$$ \Vert T_{N}( \varphi_{\alpha})\Vert \sim N^{-2\alpha}
C(\alpha)\Vert K_{\alpha}\Vert c_{1}(1),$$ 
avec $C(\alpha) =\Gamma (1+2\alpha) \frac{\sin (-\pi \alpha)}{\pi}$ et o\`u $ K_{\alpha}$ est l'op\'erateur int\'egral sur $L^2 (0,1)$ de noyau \\ $\vert x-y\vert ^{-2\alpha-1}$. Dans 
\cite{Ramb10} on utilise une ide du m\^{e}me type 
pour obtenir la valeur propre maximale d'un produit 
$T_N(\varphi_{\alpha_1}) T_N (\varphi_{\alpha_2})$
avec $-\frac{1}{2}<\alpha_1, \alpha_2<0.$ Dans ce 
cas on obtient 
$$\Vert T_N(\varphi_{\alpha_1}) T_N (\varphi_{\alpha_2})\Vert \sim N^{-2\alpha_1 -2\alpha_2}
C(\alpha_1) C(\alpha_2) c_1(1) c_2 (1) 
\Vert K_{\alpha_1, \alpha_2} \Vert$$
o $K_{\alpha_1,\alpha_2}$ est l'oprateur 
intgral sur $L^2(0,1)$ de noyau 
$$ (x,y) \rightarrow \int_0^1 \vert x-t\vert 
^{2\alpha_1 -1} \vert t-y\vert ^{2\alpha_2-1} dt.$$
Lorsque $\alpha$ est positif on cherche un quivalent de la plus grande valeur propre $ \Lambda_{N,\alpha}$ de la matrice $\left( T_N (\varphi_\alpha)\right)^{-1}$ dans le but d'obtenir la plus petite 
valeur propre de $T_N(\varphi_\alpha)$. Pour ce faire on utilise 
la fonction dfinie sur $[0,1]\times [0,1]$ si 
$\alpha>\frac{1}{2}$ et sur $[0,1]\times [0,1]\setminus 
\{(x,x)/x\in [0,1]\}$ pour $ \alpha\le \frac{1}{2}$ par
\begin{align*} G_{\alpha}(x,y) &= \frac{1}{\Gamma^2(\alpha)} x^\alpha y^\alpha \int _{\max (x,y)}^1 
\frac{ (t-x)^{\alpha-1} (t-y)^{\alpha-1} }{t^{2\alpha}} dt\quad \mathrm {si} \quad (x,y)\not=(0,0),\\
G_{\alpha}(0,0) &=0.
\end{align*}
Le lien entre cette fonction et $T_N^{-1} (\varphi_\alpha)$ a t tabli dans
 \cite{RS04}, \cite{RS1111} et 
\cite{Bot}. 
Dans \cite{BoW} B\"{o}ttcher et Widom obtiennent un \'equivalent de la valeur propre minimale  $\lambda_{N,\alpha}$ de $T_{N} (\varphi_{\alpha})$ avec $\alpha\in \mathbb N^*$. Ils utilisent le rsultat suivant, tabli dans  \cite{RS04} (voir galement \cite{Bot})

Si $\alpha\in \mathbb N^*$ on a pour tous r\'eels $x$ et $y$ avec $0\le x,y \le 1$   
$$T_{N}^{-1} (\varphi_{\alpha})_{[Nx]+1, [Ny]+1} = \frac{N^{2\alpha-1}} {c_{1}(1)} G_{\alpha} (x,y) 
+ o(N^{2\alpha-1})$$ uniform\'ement pour $x$ et $y$ dans $[0,1]$.

L'uniformit\'e de l'approximation sur $[0,1]\times [0,1]$ permet d'approcher la 
plus grande valeur propre de l'inverse par la norme d'un op\'erateur int\'egral sur $L^2(0,1)$ de noyau $G_{\alpha}(x,y)$ (cet op\'erateur \'etant de Hilbert-Schmidt et positif sa norme est aussi sa plus grande valeur propre).
B\"{o}ttcher et Widom obtiennent 
$$ \Lambda_{N,\alpha} \sim \Vert \tilde G_{\alpha}\Vert N^{2\alpha} \frac{1}{c_{1}(1)}$$
o\`u $\tilde G_{\alpha}$ est l'op\'erateur de noyau $G_{\alpha}$. \\
Dans \cite{RS1111} nous traitons ce dernier probl\`eme dans le cas o\`u $\alpha$ est un r\'eel non entier strictement sup\'erieur \`a $\frac{1}{2}$.
Quand l'exposant $\alpha$ est suprieur  $\frac{1}{2}$ mais non entier le fait que l'on n'a pas l'uniformit\'e de l'approximation sur tout $[0,1]^2$ 
(voir \cite{RS10})
emp\^{e}che d'utiliser directement l'approximation par un op\'erateur. On utilise alors une m\'ethode matricielle ( voir \cite{RS1111}) qui consiste \`a obtenir un \'equivalent quand $s$ tend vers l'infini de 
$\left(\trace \left(T_{N}(\varphi_{\alpha}) \right)^s\right)^{1/s}$ qui, puisque les valeurs propres sont positives, 
 est \'equivalent \`a $\Lambda_{N,\alpha} $. On exprime $\trace \left(T_{N}(\varphi_{\alpha}) \right)^s$ 
 au moyen des puissances de convolution $*^s G_{\alpha}(x,y)$ o\`u pour une fonction $f$ d\'efinie sur
$[0,1]^2$ on a 
\begin{align*}
\star^s f(x,y) &= \int_{0}^1\int_{0}^1 f(x,x_{1}) \int_{0}^1f(x_{1},x_{2}) \cdots\\
&\int_{0}^1 f(x_{s-1},x_{s}) f(x_{s},y) dx_{s}dx_{s-1}\cdots dx_{2}dx_{1}.
\end{align*}
On obtient alors
 \begin {equation}\label{maximal}
  \Lambda_{N,\alpha}\sim N^{2\alpha} \frac{1}{c_{1}(1)} 
 \lim_{s \rightarrow + \infty }\left(\int_{0}^1 \star^s G_{\alpha}(t,t) dt \right)^{1/s}.
 \end{equation}
 Ici nous donnons un \'equivalent de la plus petite valeur propre de $T_{N} (\varphi_{\alpha})$ avec 
$0<\alpha<\frac{1}{2}$.
Pour ce faire nous utilisons le th\'eor\`eme (voir \cite{RS10}) 
\begin{theorem} \label{Noyau}
 Pour $0<\alpha<\frac{1}{2}$ et $c_{1} \in A(\mathbb T,\frac{3}{2})$ nous avons 
$$ c_{1}(1) \left( T^{-1}_{N} (\varphi_{\alpha}) \right)_{[Nx]+1, [Ny]+1}
= N^{2\alpha-1} G_{\alpha}(x,y) +o(N^{2\alpha-1} )$$
uniform\'ement en $(x,y)$ pour $0<\delta _{1}\le x \not= y \le \delta _{2}<1$,
\end{theorem}
On rappelle que si $r \ge 1$ on pose $A(\mathbb T,r) = \{h \in L^2 (\mathbb T) \quad \mathrm{tels} \quad 
\mathrm{que} \quad
\displaystyle{\sum_{n\in \mathbb Z}} \vert \vert n+1\vert^s \hat g (n)\vert  < \infty\}$ . On peut consulter \cite {Ka} pour les espaces $A(\mathbb T,r)$.\\
Dans le cas qui nous intresse ici les indices $k,l$ pour lesquels $\left( T^{-1}_{N} (\varphi_{\alpha}) \right)_{k+1, l+1}$ ne peuvent pas 
\^{e}tre obtenus en fonction de $G_{\alpha}$ ne sont pas ngligeables dans le calcul de la trace, 
comme c'est le cas pour $\alpha>\frac{1}{2}$. Pour obtenir le th\'eor\`eme \ref{Principal} on 
dfinit la matrice $G_{N,\alpha}$ est dfinie par 
$(G_{N,\alpha})_{k+1,l+1}= N^{2\alpha-1}
\frac{1}{c_1(1)} G_{\alpha} (\frac{k}{N}, \frac{l}{N})$
si $k \not=l$ et $(G_{N,\alpha})_{k+1,k+1} =0$ et on 
montre   
$$ \Vert T_N^{-1} (\varphi_\alpha ) -G_{N,\alpha}\Vert = o(N^{2\alpha})$$
 (c'est le lemme \ref{Fonda}).
 On relie ensuite, pour terminer la dmonstration du thorme \ref{Principal}, la matrice $G_{N,\alpha}$
  l'oprateur de noyau $G_\alpha$.
Pour obtenir le lemme \ref{Fonda} nous avons besoin d'une connaissance fine des coefficients
de la diagonale de la matrice $T_N(\varphi_\alpha)$.
 Cette tude est rendue possible par une nouvelle \'ecriture des coefficients du polyn\^{o}me pr\'edicteur (on rappelle un peu plus loin le lien entre le polyn\^{o}me prdicteur de $\varphi_\alpha$ et les coefficients de la matrice $T_N^{-1} (\varphi_\alpha )$) obtenue dans \cite{RS102} \`a partir des r\'esultats de 
 \cite{RS09}. Cette criture permet d'obtenir le th\'eor\`eme \ref{Inverse1} qui relie les coefficients de $T_{N}^{-1} (\varphi_{\alpha})$ et ceux de $T_{N} (\varphi_{\alpha}^{-1})$. En gros nous obtenons 
\begin{equation} \label{morphos}
\forall k,l,\quad 0\le k,l\le N\quad \left( T_{N}^{-1} (\varphi_{\alpha}) \right)_{k+1,l+1}- \left(T_{N} (\varphi_{\alpha}^{-1})\right)_{k+1,l+1} = O(N^{2\alpha-1}).
 \end{equation}
 Ce th\'eor\`eme est \`a rapprocher de l'\'enonc\'e obtenu dans \cite{RS09} qui dit (d'une mani\`ere plus pr\'ecise qu'ici) que 
$$ \trace \left(T_{N}^{-1} (\varphi_{\alpha}) \right) - \trace \left( T_{N} (\varphi_{\alpha}^{-1})\right)
= O(N^{2\alpha}).$$
L'\'equation (\ref{morphos}) indique aussi que si le couple $(k,l)$ est au voisinage de la diagonale $T_{N}^{-1}(\varphi_{\alpha})_{k+1,l+1}$ et 
$T_{N}(\varphi^{-1}_{\alpha})_{k+1,l+1}$ ne sont pas tr\`es diff\'erents, ce qui est conforme \`a l'approximation de 
Whittle (\cite{WHI}). Pour ce qui concerne la d\'emonstration du th\'eor\`eme  \ref{Principal} l'\'egalit\'e (\ref{morphos}) intervient pour d\'emontrer le lemme \ref{Fonda}. Elle permet d'\'evaluer la diff\'erence entre les coefficients 
$\left(T_N^{-1} (\varphi_\alpha ) \right)_{k+1,l+1}$ et $(G_{N,\alpha})_{k+1,l+1}$ 
quand $\frac{\vert l- k\vert}{N} \rightarrow 0$ (dans ce cas 
le th\'eor\`eme \ref{Noyau} ne peut pas s'appliquer).\\
D'autre part la proximit entre les matrices 
$T_N^{-1}(\varphi_\alpha)$ et $G_{N,\alpha}$ permet d'approcher le produit 
$T_N^{-1}(\varphi_{\alpha_1})T_N ^{-1}(\varphi_{\alpha_2})$ par 
$G_{N,\alpha_1}G_{N,\alpha_2}$. Si 
$\alpha_1+\alpha_2>\frac{1}{2}$ on peut alors relier 
la norme de ce produit de matrices avec l'oprateur sur $L^2(0,1)$ de noyau 
$G_{\alpha_1}\star G_{\alpha_2}$.

Rappelons maintenant quelques r\'esultats et notations que nous utiliserons dans la suite de ce travail.
On sait (voir \cite{GS}) que si $c_{1}>0$ et 
$ \ln c_1$ intgrable sur $\mathbb T$ alors 
 il existe deux fonctions 
$g_{\alpha}$ et $g_{1}$ telles que 
$\varphi_{\alpha} =g_{\alpha} \overline{g_{\alpha}}$ et $c_1 =g_1 \overline{g_1}$ avec
$g_{\alpha}=(1-\chi)^\alpha g_{1}$ et $g_{1}\in H^{2+}$ (et donc $ g_{\alpha}\in H_{2+}$).
On notera ici $\beta_{u}^{(\alpha)}$ le coefficient de Fourier d'ordre $u$ de $g_\alpha^{-1}$.
On sait que si $u$ assez grand on a 
$ \beta_{u}^{(\alpha)} = \left(g_{1}(1)\right)^{-1}\frac{u^{\alpha-1} }{\Gamma (\alpha)} +o(u^{\alpha-1})$
(on pourra se r\'ef\'erer \`a \cite{Zyg2}).
Dans la suite de l'article on supposera $\beta_{0}^{(\alpha)} =1$, ce qui simplifie les notations et ne restreint pas la g\'en\'eralit\'e des r\'esultats. 
Enfin il faut remarquer que si $f$ est une fonction  valeurs relles 
$T_N (\varphi_\alpha)_{k+1,l+1} =  
\overline{ T_N (\varphi_\alpha)_{l+1,k+1}}$ et 
$T_N (\varphi_\alpha)_{k+1,l+1} =  
T_N (\varphi_\alpha)_{N-l+1,N-k+1}.$ Ces relations traduisent deux symtries sur la matrice 
$T_N (\varphi_\alpha)$ et donc sur son inverse 
$T_N^{-1} (\varphi_\alpha)$.
  Nous utiliserons enfin d'une mani\`ere d\'eterminante les propri\'et\'es des polyn\^{o}mes pr\'edicteurs. Rappelons que si $h \in L^1(\mathbb T)$
  le polyn\^{o}me pr\'edicteur de $h$ est le polyn\^{o}me trigonom\'etrique 
  dont les coefficients 
 sont obtenus en divisant les termes de la premi\`ere colonne de l'inverse de $ T_{N}(h)$ par $\left(T^{-1}_{N}(h) \right)_{1,1}^{1/2}$
(voir \cite{Ld}).\\
Rappelons ici la proprit fondamentale des polyn\^{o}mes prdicteurs ainsi que la formule de Gohberg-Semencul \cite{GoSe}.
\begin{prop}\label{fondamentale}
Si $P_N$ dsigne le polyn\^ome prdicteur de degr\'e $N$ du symbole $h$ alors 
$$
 \forall s\quad  \mathrm{tel} \quad \mathrm{que} \quad -N\le s \le N \quad \widehat {h}(s) = \widehat{\left(\frac{1}{\vert P_N \vert^2}\right)} (s).
$$ 
On a alors  
\begin{equation}
T_N(h) = T_N \left(\frac{1}{\vert P_N \vert^2}  \right).
\end{equation}
D'autre part si $Q_{N}$ est le polyn\^{o}me orthogonal associ\'e au poids  $h$ rappelons que 
\begin{equation} \label{polyortho}
Q_{N} (z) = z^N \overline{P_{N}} \left ( \frac{1}{z} \right). 
\end{equation}
 Le calcul des coefficients 
$\left( T_{N} (f)\right) ^{-1} _{k+1,l+1}$ $0\le l,k \le N$ donne donc \'egalement les coefficients des polyn\^{o}mes orthogonaux. 
\end{prop}
Rappelons la propri\'et\'e suivante 
\begin{prop}
Quelque soit l'entier naturel $N$ et le complexe $z$ appartenant \`a $\mathbb T$ on a 
$ P_{N}(z) \not=0$ et $Q_{N}(z) \not =0$.
\end{prop}
Dans la suite de ce travail nous noterons par $\chi$ la fonction $\theta \rightarrow e^{i\theta}$.
\begin{prop} (Gohberg-Semencul)  \label{THGOBSEM}
Si $K_{N} =\displaystyle{\sum_{u=0}^{N} \omega_u \chi ^u }$
un polyn\^ome trigonomtrique de degr infrieur ou gal  $N$ ne s'annulant pas sur le tore,
on a, si $0 \le k\le l \le N$ 
$$
T_N\left(\frac{1}{\vert K_{N}\vert ^2} \right)^{-1}_{k+1,l+1} 
=  \sum_{u=0}^k \bar \omega_{k-u} \omega_{l-u} -\sum_{v=0}^{k} \omega_{v+N-l} \bar \omega_{v+N-k}. 
$$
\end{prop}
 Si $f\in L^1 (\mathbb T)$ on remarque que la formule de Gohberg-Semencul et la propri\'et\'e \ref{fondamentale}  
 permettent de calculer, en toute g\'en\'eralit\'e, les coefficients $\left(T_{N}(f)\right)^{-1}_{h+1, l+1},$ $ 0 \le h\le N, \quad 0 \le l \le N$ 
 quand on conna\^{i}t les 
 coefficients $\left(T_{N}(f)\right)^{-1}_{k+1, 1} $ $0 \le k \le N.$

\section{Principaux r\'esultats}
\begin{theorem}\label{Inverse1}
Soit $\varphi_{\alpha}(e^{i \theta})=\vert 1- e^{i \theta} \vert ^{2\alpha} c_{1}(e^{i \theta})$ avec 
$0<\alpha<\frac{1}{2}$ et $c_1 \in A(\mathbb T, 
\frac{3}{2})$.
Alors il existe une fonction $h_{\alpha}$  d\'efinie  sur  
$ ]0,1]^2$  telles que pour tout r\'eel $x,y$, $0<x,y<1$ on ait,
uniform\'ement sur tout compact de $]0,1[^2$  
$$ T_{N}^{-1}(\varphi_{\alpha}) _{[Ny]+1, [Nx]+1} = \widehat {\varphi_{\alpha}^{-1}} (\vert [Ny]- [Nx] \vert)
+ \frac{N^{2\alpha-1}}{c_{1}(1) \Gamma^2(1)}
 \left( h_{\alpha} (x,y)\right) +o(N^{2\alpha-1})$$
avec 
$$ h_{\alpha}(x,y) =  h_{1,\alpha} (x,y) + h_{2,\alpha} (x,y). $$
Les fonctions $h_{1,\alpha}$ et $h_{2,\alpha}$ \'etant d\'efinies par 
$$ h_{1,\alpha} (x,y) = \int_{\min(x,y)}^{+\infty} t^{\alpha-1} (y-x+t)^{\alpha-1} dt $$
et 
\begin{align*}
 h_{2,\alpha} (x,y) &= \int_{0}^{\min(x,y)} t^{\alpha-1} (y-x+t)^{\alpha-1} \left( (1-t)^\alpha-1\right) dt 
\\
&+ \int_{0}^{\min(x,y)} (1-t)^{\alpha-1} t^\alpha (1-t-y+x)^{\alpha-1} (y-x+t)^\alpha dt 
\end{align*}
\end{theorem}

\begin{theorem} \label{Principal}
Soit $\varphi_{\alpha}(e^{i \theta})=\vert 1- e^{i \theta} \vert ^{2\alpha} c_{1}(e^{i \theta})$ avec 
$0<\alpha<\frac{1}{2}$ et 
$c_1 \in A(\mathbb T, \frac{3}{2})$. Alors si $\lambda_{min, \alpha,N}$ est  la valeur propre minimale de 
$T_N (\varphi_{\alpha})$ on a 
$$ \lambda_{min, \alpha,N} = N^{-2\alpha} \Vert  \tilde G_{\alpha} \Vert ^{-1} c_1(1) +o(N^{-2\alpha})$$
o\`u $ \tilde G_{\alpha}$ est l'op\'erateur sur $L^2 (0,1)$ de noyau $G_{\alpha}$.
\end{theorem}
Si la fonction $c_{1}$ v\'erifie les m\^{e}mes hypoth\`eses que dans le th\'eor\`eme pr\'ec\'edent nous avons \'enonc\'e dans \cite{RS10} le lemme 
\begin{lemme} \label{approx1/2} 
Si $0<\alpha<\frac{1}{2}$ et si $\frac{1}{2}-\alpha$ est suffisamment petit nous avons, avec les m\^{e}mes notations que ci-dessus 
$$ \Vert  T_N (\varphi_{\alpha}) - T_N (\varphi_{1/2}) \Vert \le K \left (\frac{1}{2} -\alpha \right) 
\Bigl \vert \ln \left (\frac{1}{2} -\alpha \right) \Bigr \vert N.$$
\end{lemme}
Ce lemme et le th\'eor\`eme \ref{Principal} permettent imm\'ediatement d'obtenir le corollaire suivant  
\begin{corollaire} 
si $\lambda_{min, 1/2,N}$ est  la valeur propre minimale de 
$T_N (\varphi_{1/2})$ on a 
$$ \lambda_{min,1/2,N} = \frac{1}{N} \Vert  \tilde G_{1/2} \Vert^{-1} c_1(1) +o(\frac{1}{N})$$
o\`u $ \tilde G_{1/2}$ est l'op\'erateur sur $L^2 (0,1)$ de noyau $G_{1/2}$.
\end{corollaire}
On peut alors donner les encadrements  
\begin{prop} \label{encadrement}
Pour $0<\alpha\le \frac{1}{2} $ on obtient 
$$  \frac{\Gamma^2(\alpha)\Gamma(2\alpha+4) }{6\Gamma (1+2\alpha)}
 \ge \Vert \tilde G_{\alpha}\Vert^{-1} \ge  \frac{\Gamma (1+\alpha) \Gamma (1-\alpha)}{\Gamma (1-2\alpha)}.$$
\end{prop}
Dans la suite si $h$ et $g$ sont deux fonctions d\'efinies  dans $L^1 ([0,1]^2)$ on note par
$f \star g (x,y)$ la fonction $(x,y)\rightarrow \int_{0}^1 f(x,t) g(t,y) dt.$
\begin{theorem}\label{PROD}
Soient $\alpha_{1}$ et $\alpha_{2}$ deux r\'eels dans $[0,\frac{1}{2}]$ avec $2\alpha_{1}+2\alpha_{2}-1>0$
et $c_{1}$ et $c_{2}$ sont deux fonctions r\'eguli\`eres appartenant \`a $A(\mathbb T, \frac{3}{2})$. 
On pose comme pr\'ec\'edemment 
$$ T_N (\varphi_{\alpha_{1}})=\vert 1- e^{i \theta} \vert ^{2\alpha_{1}} c_{1}(e^{i \theta})$$
$$T_N (\varphi_{\alpha_{2}})=\vert 1- e^{i \theta} \vert ^{2\alpha_{2}} c_{2}(e^{i \theta})$$
Alors si $\lambda_{\min,\alpha_{1},\alpha_{2}}$ d\'esigne la valeur propre minimale de 
$ T_N (\varphi_{\alpha_{1}})T_N (\varphi_{\alpha_{2}})$ on a 
$$ \lambda_{\min,\alpha_{1},\alpha_{2}} = N^{-2\alpha_{1}-2\alpha_{2}}
\Gamma^2 (\alpha_{1})\Gamma^2 (\alpha_{2})c_{1}(1) c_{2}(1)
 \Vert \tilde G_{\alpha_{1}}\star \tilde G_{\alpha_{2}}\Vert^{-1} +o(N^{-2\alpha_{1}-2\alpha_{2}}).$$
 o\`u $\tilde G_{\alpha_{1}}\star \tilde G_{\alpha_{2}}$ est l'op\'erateur de $L^2 (0,1)$ de noyau 
 $G_{\alpha_{1}}\star G_{\alpha_{2}}$
\end{theorem}
\begin{remarque} Compte tenu des rsultats 
obtenus ici et dans \cite{RS1111} on peut raisonnablement conjecturer 
qu'on peut obtenir ce rsultat pour tout couple 
d'exposants rel $\alpha_{1}>0$ et $\alpha_{2}>0$, avec $2\alpha_1+2\alpha_2-1>0$.
\end{remarque}
\begin{prop} \label{encadrement2}
Avec les m\^{e}mes hypoth\`eses que pour le th\'eor\`eme \ref{PROD} on obtient 
$$
\lambda_{\min,\alpha_{1},\alpha_{2}} = N^{-2\alpha_{1}-2\alpha_{2}}
 c_{1}(1) c_{2}(1) c_{\alpha_1,\alpha_2}
 +o(N^{-2\alpha_{1}-2\alpha_{2}})$$
 avec 
 $$\frac{\Gamma(1-\alpha_1) \Gamma(\alpha_1+1)}{\Gamma(1-2\alpha_1) }
 \frac{\Gamma(1-\alpha_2) \Gamma(\alpha_2+1)}{\Gamma(1-2\alpha_2) }
 \frac{ \min (\alpha_1,\alpha_2)}{\alpha_1+\alpha_2}
 \le c_{\alpha_1,\alpha_2} 
 $$
 et
 $$
c_{\alpha_1,\alpha_2}  \le \frac{1}{\int_0^1 (1-t)^2 t^{\alpha_1+\alpha_2}dt}
 \frac{\Gamma^2(\alpha_1)}{\int_0^1 (1-t)^2 t^{2\alpha_1}dt}
 \frac{\Gamma^2(\alpha_2)}{\int_0^1 (1-t)^2 t^{2\alpha_2}dt}.
 $$
 ou encore 
 $$
 c_{\alpha_1,\alpha_2}  \le \frac{1}{6^3}\frac{\Gamma (\alpha_{1}+\alpha_{2}+4)}{ \Gamma(\alpha_{1}+\alpha_{2}+1)}
  \frac{\Gamma(2\alpha_{1}+4) \Gamma^2 (\alpha_1)}{\Gamma (2\alpha_{1}+1)} 
  \frac{\Gamma(2\alpha_{2}+4)\Gamma^2 (\alpha_2) }{\Gamma (2\alpha_{2}+1)}.
 $$
 \end{prop}
\section{Dmonstration du th\'eor\`eme  (\ref{Inverse1})}
Dans la d\'emonstration nous allons supposer que $x<y$. Nous noterons par  
$P_{N,\alpha}=\displaystyle{\sum_{u=0}^N \gamma_{u,N}^{(\alpha)} \chi^u}$ le polyn\^{o}me pr\'edicteur de degr\'e $N$ de la fonction $\varphi_{\alpha}$. Pour en savoir plus sur les polyn\^{o}mes pr\'edicteurs on pourra se r\'ef\'erer \`a \cite{Ld} ou \cite{RS04}.
Nous utiliserons  le r\'esultat suivant \'etabli dans \cite{RS102} 
\begin{theorem} \label{theoremepolypredi}
On consid\`ere une fonction $\varphi_{\alpha}$ v\'erifiant les hypoth\`eses du th\'eor\`eme 
\ref{Inverse1}.
Alors il existe un entier $n_{1}$, ind\'ependant de $N$, tel que 
$$ \gamma_{k,N}^{(\alpha)} = \beta_{k}^{(\alpha)} (1-\frac{k}{N})^\alpha \left(1+o(1)\right)$$
pour tout entier $k \in [0, N -n_{1} ]$, uniform\'ement par rapport \`a $N$.
\end{theorem}
\begin{remarque}\label{Rem2}
Dans la pratique $n_{1}$ est choisi par rapport \`a un r\'eel $\epsilon>0$ de mani\`ere \`a ce que pour tout entier $0\le u\ge N-n_{1}$ on ait 
$\beta_{u}^{(\alpha)} = \frac{1}{g_{1}(1)} \frac{u^{\alpha-1}}{\Gamma(\alpha)} (1+r_{u})$ avec la pr\'ecision $\vert r_{u}\vert < \epsilon.$
\end{remarque}
\begin{remarque}
Ce th\'eor\`eme peut alors se lire, $n_{1}$ \'etant comme dans la remarque \ref{Rem2} \\
$ \forall \epsilon >0 \quad  \exists N_{0} \quad \mathrm{t. q.}\quad  \forall N \ge N_{0}
\quad \forall k, \quad 0\le k \le N-n_{1} \quad  \exists R_{k}, \quad \vert R_{k} \vert \le \varepsilon$
 tel que  
$$ \gamma_{k,N}^{(\alpha)} = \beta_{k}^{(\alpha)} (1-\frac{k}{N})^\alpha \left(1+R_{k}\right).$$
\end{remarque}
Nous allons aussi utiliser le th\'eor\`eme suivant qui nous permet de pr\'eciser les coefficients 
$\gamma_{k,N}^{(\alpha)}$ quand $N \rightarrow +\infty$ (voir \cite{RS1111}).
\begin{theorem} \label{rappel}
Soit $\varphi_{\alpha}$ est une fonction v\'erifiant les hypoth\`eses du th\'eor\`eme et telle que 
$\beta_{0}^{(\alpha)} =1$. Si $k$ un entier 
tel que $\displaystyle{\frac{k}{N} \rightarrow 0}$ on a  
$$ \gamma_{N-k}^{(\alpha)} =  \beta_{k}^{(\alpha+1)} \frac{\alpha}{N} \left( 1+o(1)\right),$$ 
o\`u $ \beta_{k}^{(\alpha+1)} $ est le coefficient de Fourier d'ordre $k$ de la 
fonction $\varphi_{\alpha+1}=\vert 1-\chi\vert ^{2(\alpha+1)} c_{1}.$
\end{theorem}
\begin{remarque} On peut remarquer que si 
$\frac{k}{N} \rightarrow 0$ avec $k$ suprieur au $n_1$ de la remarque 1 les thormes \ref{theoremepolypredi} et \ref{rappel} sont compatibles pour calculer 
$\gamma_{N+1-k}^{(\alpha)} $.
\end{remarque}
Nous allons utiliser maintenant la formule de Gohberg-Semencul (voir l'intoduction), qui nous permettra de calculer les coefficients 
de la matrice $T_{N}^{-1}(\varphi_{\alpha})$ en fonction des coefficients $\gamma_u^{(\alpha)}$.
On a, en posant $k=[Nx]$, $l=[Ny]$ et en supposant $k\le l$
\begin{equation} \label{GOBSEM}
\left( T_{N}^{-1}(\varphi_{\alpha})\right) _{k+1,l+1} = 
 \sum_{u=0} ^k \overline{\gamma _{k-u,N}^{(\alpha)}}
\gamma _{l-u,N}^{(\alpha)}-\sum_{v=0}^{k}
\gamma_{v+N-l,N}^{(\alpha)}  \overline{\gamma _{v+N-k,N}^{(\alpha)}}.
\end{equation}
Il vient alors 
\begin{align*}
\sum_{u=0}^k \overline{\gamma_{u,N}^{(\alpha)}} \gamma_{l-k+u,N}^{(\alpha)}
&= \sum_{u=0}^k \overline{\beta_{u}^{(\alpha)}} \beta_{l-k+u}^{(\alpha)} + \sum_{u=0}^k \overline{\beta_{u}^{(\alpha)}} \left(\gamma_{l-k+u,N}^{(\alpha)}- 
\beta_{l-k+u}^{(\alpha)}\right)+\\
& +\sum_{u=0}^k \left( \overline{\gamma_{u,N}^{(\alpha)} - \beta_{u}^{(\alpha)}}\right) \gamma_{l-k+u,N}^{(\alpha)}.
\end{align*}
Nous pouvons \'ecrire, si $k$ assez grand pour que $\beta_{u}^{(\alpha)}$ puisse \^{e}tre remplac\'e par son asymptotique pour $u \ge k$ 
\begin{align*} 
 \sum_{u=0}^k \overline{\beta_{u}^{(\alpha)}} \beta_{l-k+u}^{(\alpha)} &=
  \sum_{u=0}^{+ \infty} \overline{\beta_{u}^{(\alpha)}} \beta_{l-k+u}^{(\alpha)}
-  \sum_{u=k+1}^{+ \infty} \overline{\beta_{u}^{(\alpha)}} \beta_{l-k+u}^{(\alpha)}\\
&= 
\widehat{\varphi_{\alpha}^{-1}} (l-k) - \frac{N^{2\alpha-1} }{\Gamma^2(\alpha) c_{1}(1)} 
I_{1,\alpha} (x,y) + o(N^{2\alpha-1})
\end{align*}
avec 
$$
I_{1,\alpha}(x,y) = \int_{x}^{+ \infty} t^{\alpha-1} (y-x+t)^{\alpha-1} dt.$$
Nousavonsensuite 
 $$ \sum_{u=0}^k \overline{\beta_{u}^{(\alpha)}} \left(\gamma_{l-k+u,N}^{(\alpha)}- 
\beta_{l-k+u}^{(\alpha)}\right) = \frac{N^{2\alpha-1} }{\Gamma^2(\alpha) c_{1}(1)} 
I_{2,\alpha} (x,y) + o(N^{2\alpha-1})$$
avec 
$$ I_{2,\alpha} (x,y) = \int_{0}^{x} t^{\alpha-1} (y-x+t)^{\alpha-1} \left( (1-y+x-t)^\alpha-1\right) dt.$$
En effet soti $k_{0}$ un entier ind\'ependant de $N$ tel que pour tout $u\ge k_{0}$ l'on puisse remplacer 
$\beta_{u}^{(\alpha)}$ par son asymptotique. 
On a $k_{0}\le k$ si $x>0$ et $N$ assez grand. On peut alors \'ecrire 
\begin{itemize}
\item [$i)$]
si $l-k>k_{0}$
\begin{align*}
 \sum_{u=0}^{k_{0}} \overline{\beta_{u}^{(\alpha)}} \left(\gamma_{l-k+u,N}^{(\alpha)}- 
\beta_{l-k+u}^{(\alpha)}\right) &=  \sum_{u=0}^{k_{0}} \overline{\beta_{u}^{(\alpha)}} 
\beta_{l-k+u}^{(\alpha)} \left( 1-(1-\frac{l-k+u}{N})^\alpha\right)\\
&\sim N^{\alpha-1} \frac{(y-x)^{\alpha-1}}{\Gamma (\alpha) g_{1}(1)} \left(1-(1-y+x)^\alpha\right)
\sum_{u=0}^{k_{0}} \overline{\beta_{u}^{(\alpha)}} \\
&=O(N^{\alpha-1}) =o(N^{2\alpha-1}).
\end{align*}
\item [$ii)$]
Si $l-k \le k_{0}$ alors 
$$ \sum_{u=0}^{k_{0}} \overline{\beta_{u}^{(\alpha)}} \left(\gamma_{l-k+u,N}^{(\alpha)}- 
\beta_{l-k+u}^{(\alpha)}\right) = O\left( \frac{1}{N}\right) =o(N^{2\alpha-1}).$$
\end{itemize}
Etnfin on obtient 
$$  \sum_{k_{0}}^k \overline{\beta_{u}^{(\alpha)}} \left(\gamma_{l-k+u,N}^{(\alpha)}- 
\beta_{l-k+u}^{(\alpha)}\right)=\frac{N^{2\alpha-1} }{\Gamma^2(\alpha) c_{1}(1)} 
I_{2,\alpha} (x,y) + o(N^{2\alpha-1})$$
avec la formule d'Euler et Mac-Laurin.
Ces m\^{e}mes m\'ethodes nous donnent d'une part  
$$ \sum_{u=0}^k \overline{\gamma_{u,N}^{(\alpha)} - \beta_{u}^{(\alpha)}} \gamma_{l-k+u,N}^{(\alpha)} =
 \frac{N^{2\alpha-1} }{\Gamma^2(\alpha) c_{1}(1)} 
I_{3,\alpha} (x,y) + o(N^{2\alpha-1})$$
avec 
$$ I_{3,\alpha} (x,y) = \int_{0}^{x} (y-x+t)^{\alpha-1} (1-y+x-t)^\alpha t^{\alpha-1}  \left( (1-t)^\alpha-1\right) dt,
$$
et d'autre part 
$$ \sum_{u=0}^k  \gamma_{N-l+u,N}^{(\alpha)}\overline{\gamma_{N-k+u,N}^{(\alpha)}}
= \frac{N^{2\alpha-1} }{\Gamma^2(\alpha) c_{1}(1)} 
I_{4,\alpha} (x,y) + o(N^{2\alpha-1})$$
avec 
$$ I_{4,\alpha} (x,y) = \int_{0}^{x} t^{\alpha-1} (1-t)^\alpha (y-x+t)^{\alpha-1}  (1-y+x-t)^\alpha dt. $$
On peut remarquer que l'uniformit\'e annonc\'ee est fournie par l'uniformit\'e de l'approximation de 
$\beta_{u}^{(\alpha)} $ par $\frac{u^{\alpha-1}}{\Gamma(\alpha) g_{1}(1)}$ et par le reste de 
la formule d'Euler et Mac-Laurin.
Ceci ach\`eve de prouver le th\'eor\`eme \ref{Inverse1}.
Pour obtenir le th\'eor\`eme \ref{Principal} nous avons besoin d'une \'etude plus fine de 
certains  \'el\'ements de la matrice 
$T_{N}^{-1}(\varphi_{\alpha}).$ Cela va \^{e}tre le but du th\'eor\`eme \ref{Inverse2} que nous allons
\'enoncer, puis d\'emontrer.
\begin{theorem}\label{Inverse2} 
Soit un r\'eel $\epsilon$ strictement positif. Si $n_{\epsilon}$ est un entier naturel tel que 
$$ \forall u \ge n_{\epsilon} \quad \beta_{u}^{(\alpha)} = \frac{u^{\alpha-1}}{\Gamma (\alpha) g_{1}(1)} 
\left( 1+ R(u)\right) \quad \mathrm{avec} \quad \vert R(u)\vert <\epsilon.$$
Alors pour tout r\'eel $\delta \rightarrow 0$ avec $N\delta >n_{\epsilon}$ il existe une constante 
$C_{1,\alpha}$ qui ne d\'epend que de $\alpha$ telle que
$$ \left( T_{N} ^{-1}(\varphi_{\alpha}) \right)_{k+1,l+1} \le C_{1,\alpha} \vert l-k\vert^{\alpha-1} \left( N\delta \right)^\alpha
$$
pour tout couple d'entiers naturels $k,l$ avec 
$  0 \le \min (k,l) <N\delta $ et  $2 N \delta < \max (k,l) <  N- 2 N \delta $.
\end{theorem}
\begin{remarque}
La conclusion du th\'eor\`eme \ref{Inverse2}  peut \'egalement s'\'enoncer \\
il existe une constante 
$C_{1,\alpha}$ qui ne d\'epend que de $\alpha$ telle que
$$ \left( T_{N} ^{-1} (\varphi_{\alpha}) \right)_{k+1,l+1} \le C_{2,\alpha} \vert \widehat{\varphi_{\alpha/2}} (k-l)\vert \left( N\delta \right)^\alpha
$$
pour tout couple d'entiers naturels $k,l$ avec 
$  0 \le \min (k,l) <N\delta $ et  $2 N \delta < \max (k,l) <  N- 2 N \delta $.
\end{remarque}
\section{Dmonstration du thorme \ref{Inverse2}}
On reprend la formule de Gohberg-Semencul (formule \ref{GOBSEM}) avec $\min (k,l) =k$ et 
$\max (k,l) =l$. Les autres cas se d\'eduisent de ce cas l\`a en utilisant les sym\'etries de la matrices de
Toeplitz qui se transmettent \`a son inverse. On repart de la dcomposition 
\begin{align*}
\sum_{u=0}^k \overline{\gamma_{u,N}^{(\alpha)}} \gamma_{l-k+u,N}^{(\alpha)}
&= \sum_{u=0}^k \overline{\beta_{u}^{(\alpha)}} \beta_{l-k+u}^{(\alpha)} + \sum_{u=0}^k \overline{\beta_{u}^{(\alpha)}} \left(\gamma_{l-k+u,N}^{(\alpha)}- 
\beta_{l-k+u}^{(\alpha)}\right)+\\
& +\sum_{u=0}^k \left( \overline{\gamma_{u,N}^{(\alpha)} - \beta_{u}^{(\alpha)}}\right) \gamma_{l-k+u,N}^{(\alpha)}.
\end{align*}
Ecrivons 
$$
\sum_{u=0}^k \overline{\beta_{u}^{(\alpha)}} \beta_{l-k+u}^{(\alpha)} =  \sum_{u=0}^k \overline{\beta_{u}^{(\alpha)}}\left(  \frac{(l-k+u)^{\alpha-1} }{\Gamma (\alpha)g_1(1) }\left (1+o(1)\right)\right)$$
ce qui donne la majoration 
$$ \Bigr \vert \sum_{u=0}^k \overline{\beta_{u}^{(\alpha)}} \beta_{l-k+u}^{(\alpha)} \Bigl \vert 
\le \frac{(l-k)^{\alpha-1} }{\Gamma (\alpha)
\vert g_1(1)\vert }\sum_{u=0}^{N \delta} \vert \beta_{u}^{(\alpha)}\vert$$
ou encore 
$$ \Bigr \vert \sum_{u=0}^k \overline{\beta_{u}^{(\alpha)}} \beta_{l-k+u}^{(\alpha)} \Bigl \vert 
\le \frac{(l-k)^{\alpha-1} (N\delta)^\alpha}{\Gamma^2(\alpha)\vert c_1(1)\vert }$$
ce qui est aussi 
$$ \Bigr \vert \sum_{u=0}^k \overline{\beta_{u}^{(\alpha)}} \beta_{l-k+u}^{(\alpha)} \Bigl \vert 
\le C \vert  \widehat{\varphi_{\alpha/2}^{-1}}(l-k) \vert (N\delta)^\alpha
$$
avec $C=\left( \Gamma (1-\alpha) \sin (\frac{\pi\alpha}{\pi} \right)^{-1} \frac{1}{\Gamma ^2 (\alpha) c_{1}(1)}$.
De m\^{e}me 
$$ \sum_{u=0}^k \overline{\beta_{u}^{(\alpha)}} \left(\gamma_{l-k+u,N}^{(\alpha)}- 
\beta_{l-k+u}^{(\alpha)}\right)
=  \sum_{u=0}^k \overline{\beta_{u}^{(\alpha)}} 
\left(  \frac{(l-k+u)^{\alpha-1} }{\Gamma (\alpha)g_1(1) }\left( (1-\frac{l-k+u}{N})^\alpha-1\right) \right)
\left(1+o(1)\right).$$
En procdant comme prcdemment on obtient la majoration 
$$ \Bigl \vert  \sum_{u=0}^k \overline{\beta_{u}^{(\alpha)}} \left(\gamma_{l-k+u,N}^{(\alpha)}- 
\beta_{l-k+u}^{(\alpha)}\right)\Bigr \vert 
\le 
\vert \frac{(l-k)^{\alpha-1} }{\Gamma (\alpha)g_1(1) }\vert \vert (1-\frac{l-k}{N})^\alpha-1\vert 
\Bigl \vert  \sum_{u=0}^k \overline{\beta_{u}^{(\alpha)}} \Bigr\vert.$$
Et finalement 
$$ \Bigl \vert  \sum_{u=0}^k \overline{\beta_{u}^{(\alpha)}} \left(\gamma_{l-k+u,N}^{(\alpha)}- 
\beta_{l-k+u}^{(\alpha)}\right)\Bigr \vert 
\le  \frac{(l-k)^{\alpha-1} (N\delta)^\alpha}{\Gamma
^2(\alpha)\vert c_1(1)\vert }$$
ou aussi 
$$ \Bigl \vert  \sum_{u=0}^k \overline{\beta_{u}^{(\alpha)}} \left(\gamma_{l-k+u,N}^{(\alpha)}- 
\beta_{l-k+u}^{(\alpha)}\right)\Bigr \vert 
\le C\vert  \widehat{\varphi_{\alpha/2}^{-1}}(l-k) \vert (N\delta)^\alpha.
$$
Enfin nos obtenons, avec les m\^{e}mes 
procds,
 $$\Bigl \vert
 \sum_{u=0}^k \left( \overline{\gamma_{u,N}^{(\alpha)} - \beta_{u}^{(\alpha)}}\right) \gamma_{l-k+u,N}^{(\alpha)}\Bigr \vert 
 \le 2 \frac{(l-k)^{\alpha-1} (N\delta)^\alpha}{\Gamma^2(\alpha)\vert c_1(1)\vert }$$
 ou
  $$\Bigl \vert
 \sum_{u=0}^k \left( \overline{\gamma_{u,N}^{(\alpha)} - \beta_{u}^{(\alpha)}}\right) \gamma_{l-k+u,N}^{(\alpha)}\Bigr \vert 
 \le 2
 C \vert  \widehat{\varphi_{\alpha/2}^{-1}}(l-k) \vert (N\delta)^\alpha. 
$$
Reste \`a traiter le deuxime terme de la formule de Gohberg-Semencul (formule \ref{GOBSEM}). A savoir 
$\displaystyle{ \sum_{u=0}^k \gamma_{N-k+u,N}^{(\alpha)} \gamma_{N-l+u,N}^{(\alpha)} }$. 
En utilsant le th\'eor\`eme \ref{rappel} on obtient 
\begin{align*}
\Bigl \vert \sum_{u=0}^k \gamma_{N-k+u,N}^{(\alpha)} \gamma_{N-l+u,N}^{(\alpha)}\Bigr \vert 
&\le M_{\alpha} \sum_{u=0}^k\Bigl  \vert \frac{\beta^{(\alpha+1)} _{k-u}}{N} \Bigr \vert (N-l+u)^{\alpha-1}
 \left( \frac{l-u}{N} \right)^\alpha\\
 &\le M'_{\alpha} N^{-1} (N\delta )^{\alpha-1} \le M'_{\alpha}  \vert l-k\vert ^{\alpha-1} (N\delta )^{\alpha-1}.
 \end{align*}
 En remarquant que $M_{\alpha}$ et $M'_{\alpha}$ ne d\'ependent que de $\alpha$ ceci termine la d\'emonstration du th\'eor\`eme.

 \section{D\'emonstration du th\'eor\`eme \ref{Principal}}
 \subsection {R\'esultats pr\'eliminaires}
Nous allons d'abord devoir obtenir les quatre lemmes suivants 
\begin{lemme} \label{LEMME1}
Si $0<\alpha$, $\alpha\neq1$ et $0 \le x \not=y \le 1$ on a 
$$ G_{\alpha}(x,y) \le C_{\alpha} \vert x-y \vert ^{2\alpha-1}\quad \mathrm{avec} \quad C_{\alpha}=\frac{ \Gamma (1-2\alpha) } {\Gamma (1-\alpha)\Gamma (\alpha) } $$
\end{lemme}
\begin{lemme} \label{LEMME2}
Si $0<\alpha<\frac{1}{2}$ et $0 < x \not=y \le 1$ il existe une constante $H$, ind\'ependante de $x$ et $y$ telle que
$$\vert h_{\alpha}(x,y)\vert  \le H\vert y-x\vert ^{\alpha-1}.$$
\end{lemme}

\begin{preuve}{du lemme \ref{LEMME1}} 
Pour la d\'emonstration de ce lemme nous supposerons $0\le x<y\le 1$. Il est alors clair que 
$$ \int_{y}^1 \frac{(t-x)^{\alpha-1} (t-y)^{\alpha-1} }{t^{2\alpha}} dt \le 
\frac{1}{y^{2\alpha}} \int_{y}^1 (t-x)^{\alpha-1} (t-y)^{\alpha-1} dt $$
Nous allons nous concentrer sur l'int\'egrale $ \int_{y}^1 (t-x)^{\alpha-1} (t-y)^{\alpha-1} dt.$
En utilisant des changements de variables successifs nous obtenons :
\begin{align*} 
\int_{y}^1 (t-x)^{\alpha-1} (t-y)^{\alpha-1} dt &=
(y-x)^{\alpha-1} \int_{y-x}^{1-x} h^{\alpha-1} \left( \frac{h}{y-x} -1 \right)^{\alpha-1} dh  \\
&= (y-x) ^{2\alpha-1}\int_{1}^{\frac{1-x}{y-x}} u^{\alpha-1} (u-1)^{\alpha-1} du  \\
&=  (y-x) ^{2\alpha-1}\int_{\frac{y-x}{1-x} } ^1 v^{-2\alpha} (1-v)^{\alpha-1} dv.
\end{align*}
Et puisque 
$$ \int_{0}^1v^{-2\alpha} (1-v)^{\alpha-1} dv = \frac{\Gamma (1-2\alpha) \Gamma (\alpha)}
{\Gamma (1-\alpha)}$$
nous pouvons \'ecrire 
\begin {align*}
G_{\alpha}(x,y) &\le \frac{x^\alpha}{y^\alpha}  \frac{\Gamma (1-2\alpha) \Gamma (\alpha)}
{\Gamma (1-\alpha)}(y-x) ^{2\alpha-1}\\
&\le  \frac{\Gamma (1-2\alpha) }
{\Gamma (1-\alpha)\Gamma (\alpha)}(y-x) ^{2\alpha-1}.
\end{align*}
\end{preuve}
\begin{preuve} {du lemme \ref{LEMME2}}
En remarquant que si $t>1$ alors $ y-x+t > t(y-x)$ nous avons, en supposant encore cette fois que 
$0<x< y <1$,  
$$\vert h_{1,\alpha}(x,y)\vert = \int_{x}^{+ \infty} t^{\alpha-1} (y-x+t)^{\alpha-1} dt \le \left(  \int_{0}^1 t^{\alpha-1} dt + \int_{1}^{+ \infty} 
t^{2\alpha-1}  dt \right) (y-x)^{\alpha-1} dt.$$
Nous avons d'autre part 
$$  \int_{0}^{x} t^{\alpha-1} (y-x+t)^{\alpha-1}\left( (1-t)^\alpha -1\right) dt 
\le (y-x)^{\alpha-1} \int_{0}^1 t^{\alpha-1} \left( (1-t)^\alpha-1\right) dt.$$
Ensuite, en remarquant que l'on a $t<x\Rightarrow 1-t > 1-x>y-x$ on peut \'ecrire

\begin{align*}
 &\Bigl \vert \int_{0} ^x (1-t)^{\alpha-1} t^\alpha (1-t-y+x)^{\alpha-1} (y-x+t)^\alpha dt  \Bigr \vert
\\
&\le
 (y-x)^{\alpha-1} \int_{0}^x  t^{\alpha} (1-t-y+x)^{\alpha-1} dt \\
 \le & (y-x)^{\alpha-1} \int_{0}^x   (1-t-y+x)^{\alpha-1} dt \\
 \le&   (y-x)^{\alpha-1}\left( \frac{(1-y)^\alpha}{\alpha}
 +  \frac{(1-y+x)^\alpha}{\alpha}\right)\le \frac{2}{\alpha}.
(y-x)^{\alpha-1}
\end{align*}
C'est \`a dire que $$\vert h_{2,\alpha}(x,y) \le (y-x)^{\alpha-1} \left( \frac{2}{\alpha} +   \int_{0}^1 t^{\alpha-1} \left( (1-t)^\alpha-1\right) dt\right).$$
 
Ce qui donne la majoration annonc\'ee.
\end{preuve}
\subsection{Un Lemme d'approximation.}
Soit $\delta$ un r\'eel suffisamment petit.
Pour la suite de la d\'emonstration nous allons introduire les sous-ensembles  suivants de $[0,N]^2\cap \mathbb N$,
 les intervalles utiliss ici tant des intervalles de $\mathbb N$ et $N_{1}$ d\'esignant la partie enti\`ere de 
 $N\delta $.
\begin{enumerate}
\item
$$ I_{1,\delta }=[0, 2N_{1}]^2, \quad I_{2,\delta } =[N-2 N_{1} ,N]^2,$$
$$ I_{3,\delta } =[N-2
N_{1} ,N]
\times [0, 2N_{1}], \quad I_{4,\delta}= [0,2N_{1} ]\times [N-2N_{1} ,N],$$
\item
$$ L_{1,\delta } = \{ (i,j)/ 0 \le i \le N_{1}, \,2 N_1\le j
\le N-2N_{1}\},$$
$$ L_{2,\delta } =\{ (i,j) / 2N_{1} \le i<N-2N_{1}, \, 
0\le j \le N_{1}  \},$$
$$ L_{3,\delta }= \{  (i,j) / N- N_{1} \le i \le N, \, 
2 N_{1} \le j \le N- 2N_{1}\},$$
$$ L_{4,\delta }= \{  (i,j) /  2N_{1} \le i \le N -2N_{1}, \, 
N-N_{1} \le j \le N  \}.$$
\item
$$\Delta_{\delta } =\{(i,j)\in \mathbb N^2/ 0\le \vert i-j \vert \le N_{1} \}.$$
\item
$$ D_{\delta }= \Delta_{\delta } \setminus \left(  I_{1,\delta } \cup  I_{2,\delta }\right).$$
\item
$$ C_{\delta }= [0,N]^2 \setminus J_{\delta }$$
avec $$ J_{\delta}= \Delta_{\delta }\cup \left ( \cup_{h=1}^4 I_{h,\delta }\right)
 \cup \left ( \cup_{h=1}^4 L_{h,\delta }\right).
 $$
\end{enumerate}
Le thorme \ref{Principal} est alors la consquence du lemme 
\begin{lemme} \label{Fonda}
Si $G_{N,\alpha}$ dsigne la matrice $(N+1) \times (N+1)$ dfinit pour tout entier $k,l$ $0\le k,l\le N$ par $(G_{N,\alpha}) _{k+1,l+1} =
N^{2\alpha-1}
\frac{1}{c_1(1)} G_{\alpha} (\frac{k}{N}, \frac{l}{N})$
si $k \not=l$ et $(G_{N,\alpha})_{k+1,k+1} =0$.
Alors 
$$ \Vert T_N^{-1} (\varphi_\alpha) -G_{N,\alpha}\Vert = o(N^{2\alpha}).$$
\end{lemme}
\begin{remarque} 
On sait que $T_N^{-1} (\varphi_\alpha) $ est une matrice diagonalisable, dont les valeurs propres sont srtictement positives et on a 
$\Vert T_N^{-1} (\varphi_\alpha) \Vert = \Lambda_{N,\alpha}$ si $\Lambda_{N,\alpha}$ d\'esigne la plus grande de ces valeurs propres. D'autre part $G_{N,\alpha}$  est une matrice sym\'etrique donc diagonalisable \`a valeurs propres r\'eelles.  Si $\tilde\Lambda_{N,\alpha}$ d\'esigne la valeur propre maximale de cette matrice il est clair que 
$\Vert  G_{N,\alpha} \Vert = \vert \tilde\Lambda_{N,\alpha} \vert$. Le lemme \ref{Fonda} implique donc que 
$\Bigr \vert \Lambda_{N,\alpha}  - \vert \tilde \Lambda_{N,\alpha}\vert \Bigl \vert =o(N^{2\alpha}).$ La suite de la d\'emonstration, apr\`es la preuve du lemme \ref{Fonda}, sera donc consacr\'ee \`a l'estimation de 
$\vert\tilde\Lambda_{N,\alpha}\vert$.
\end{remarque}
\begin{preuve}{du lemme \ref{Fonda}}
Posons $T_N^{-1} (\varphi_{\alpha}) -G_{N,\alpha,}
= H_{N,\alpha}$ et consid\'erons un r\'eel $\delta >0$ qui tend vers z\'ero.
On a 
$$ \Vert H_{N,\alpha}\Vert = \max_{\Vert x\Vert = \Vert y\Vert =1} \langle H_{N,\alpha} (x) \vert y \rangle.$$
Ecrivons 
$$  \langle H_{N,\alpha} (x) \vert y \rangle =
\sum_{i=0}^N \left( \sum_{j=0}^N (H_{N,\alpha})_{i+1,j+1} x_{j+1} \right) y_{i+1}.$$
Nous allons maintenant utiliser la dcomposition 
\begin{align*} 
& \sum_{i=0}^N \left( \sum_{j=0}^N (H_{N,\alpha})_{i+1,j+1} x_{j+1} \right) y_{i+1} =
\sum _{(i,j) \in C_{\delta } } 
(H_{N,\alpha})_{i+1,j+1} x_{j+1}  y_{i+1} \\
& 
+ \sum _{(i,j) \in D_{\delta } } 
(H_{N,\alpha})_{i+1,j+1} x_{j+1}  y_{i+1} 
+ \sum_{k=1}^4 \sum _{(i,j) \in I_{k,\delta } } 
(H_{N,\alpha})_{i+1,j+1} x_{j+1}  y_{i+1},\\
&
+\sum_{k=1}^2 \sum _{(i,j) \in L_{k,\delta } } 
(H_{N,\alpha})_{i+1,j+1} x_{j+1}  y_{i+1}.
\end{align*}
Avec le thorme \ref{Noyau} et la d\'efinition de $G_{N,\alpha}$ on vrifie facilement que 
$$
 \Bigl \vert \sum_{(i,j)\in C_{\delta }} (H_{N, \alpha,\delta})_{i+1,j+1} x_{i+1} y_{j +1}\Bigr \vert 
 \le   \max_{(i,j) \in C_{\delta }}\vert  (H_{N, \alpha})_{i+1,j+1} \vert \sum_{0\le i \le N, 0 \le j \le N} \vert x_{i+1} y_{j+1} \vert = o(N^{2\alpha-1}).
 $$
Nous avons d'autre part 
\begin{align*}
\Bigl \vert \sum_{(i,j)\in D_{\delta }} (H_{N, \alpha})_{i+1,j+1} x_{i+1} y_{j +1}\Bigr \vert \le &
\Bigl \vert \sum_{(i,j)\in D_{\delta }} \left( T_{N}^{-1} (\varphi_{\alpha})\right)_{i+1,j+1} x_{i+1} y_{j +1}\Bigr \vert\\
&+ \Bigl \vert \sum_{(i,j)\in D_{\delta }} (G_{N, \alpha,\delta})_{i+1,j+1} x_{i+1} y_{j +1}\Bigr \vert
\\
 \le&\sum_{(i,j)\in D_{\delta }} \Bigl \vert \widehat {\varphi_{\alpha}^{-1}}(i-j)  +N^{2\alpha-1} 
 h_{\alpha}(\frac{i+1}{N},\frac{j+1}{N}) \Bigr \vert\vert x_{i+1}\vert
 \vert  y_{j +1}\vert\\
 +& \sum_{(i,j)\in D_{\delta }}\Bigl \vert  (G_{N, \alpha,\delta})_{i+1,j+1} \Bigr \vert \vert x_{i+1} \vert \vert y_{j +1}\vert
\end{align*}

On a vu dans le lemme \ref{LEMME1} que 
$$ (G_{N, \alpha})_{i+1,j+1} \le C_{\alpha} \vert \frac {k-i}{N} \vert ^{2\alpha-1} N^{2\alpha-1}.$$
D'o\`u 
\begin{eqnarray*}
 \sum_{(i,j)\in D_{\delta }}\Bigl \vert  (G_{N, \alpha})_{i+1,j+1} \Bigr \vert \vert x_{i+1} \vert \vert y_{j +1}\vert 
&\le & C_{\alpha} N^{2\alpha-1}  
\sum_{(i,j)\in D_\delta} \Bigl \vert \frac{i-j}{N}\Bigr \vert ^{2\alpha-1} 
\vert x_{i+1} \vert \vert y_{j +1}\vert \\
&\le &
 C_{\alpha} N^{2\alpha}\int_{-\delta }^{\delta } t^{2\alpha-1} dt = O( N^{2\alpha} \delta ^{2\alpha} )=o(N^{2\alpha}).
 \end{eqnarray*}
De m\^{e}me on a,  en utilisant cette fois le lemme \ref{LEMME2}
 $$
 \sum_{(i,j)\in D_{\delta }}\Bigl \vert N^{2\alpha-1} 
 h_{\alpha}(\frac{i+1}{N},\frac{j+1}{N})   \Bigr \vert \vert x_{i+1} \vert \vert y_{j +1}\vert 
 \le C'_{\alpha} N^{2\alpha}\int_{-\delta }^{\delta } t^{\alpha-1} dt = O( N^{2\alpha} \delta ^{\alpha} )=o(N^{2\alpha}).
$$
Enfin, puisque $\widehat {\varphi_{\alpha}} (u) =O(u^{2\alpha-1})$ si $u$ assez grand on a :
$$ \Bigl \vert \sum_{(i,j)\in D_{\delta }} \Bigl \vert \widehat {\varphi_{\alpha}^{-1}}(i-j) \Bigr \vert x_{i+1} y_{j +1}\Bigr \vert
=O(N^{2\alpha} \delta ^{2\alpha}) =o(N^{2\alpha}).$$
Consid\'erons maintenant, quelque soit l'entier $k$, $1\le k\le 4,$ la quantit\'e 
 \begin{align*}
\Bigl \vert \sum_{(i,j)\in L_{k,\delta }} (H_{N, \alpha})_{i+1,j+1} x_{i+1} y_{j +1}\Bigr \vert \le &
\Bigl \vert \sum_{(i,j)\in L_{k,\delta }} \left( T_{N}^{-1} (\varphi_{\alpha})\right)_{i+1,j+1} x_{i+1} y_{j +1}\Bigr \vert\\
&+ \Bigl \vert \sum_{(i,j)\in L_{k,\delta }} (G_{N, \alpha})_{i+1,j+1} x_{i+1} y_{j +1}\Bigr \vert
\\
\end{align*}
 
 On obtient, avec le th\'eor\`eme  \ref{Inverse2} 
$$ 
\Bigl \vert \sum _{(i,j) \in L_{k,\delta} } \left(T_{N}^{-1} (\varphi_{\alpha})\right)_{i+1,j+1} x_{j+1}  y_{i+1}\Bigr \vert 
\le
C_{1,\alpha} (N\delta )^{\alpha}  \sum _{(i,j) \in L_{k,\delta } }\vert i-j\vert^{\alpha-1} \vert x_{j+1}  y_{i+1} \vert.$$
En \'ecrivant, si par exemple $k=1$,  
$$
\sum _{(i,j) \in L_{1,\delta } }\vert i-j\vert^{\alpha-1} \vert x_{j+1}  y_{i+1} \vert =
\sum_{v=N_1 }^{N-2N_1 } v^{\alpha-1} \left( \sum_{i=0 } ^{N_1 } \vert y _{v+i}\vert 
\vert x_{i}\vert \right)$$
on obtient, puisque $\Vert x\Vert =\Vert y \Vert =1$, 
$$
\Bigl \vert \sum _{(i,j) \in L_{1,\delta } } \left(T_{N}^{-1} (\varphi_{\alpha})\right)_{i+1,j+1} x_{j+1}  y_{i+1}\Bigr \vert 
= O\left((N^{2\alpha} \delta ^\alpha\right).$$

 Toujours en utilisant le lemme \ref{LEMME1} on obtient 
 $$
 \sum_{(i,j)\in L_{1,\delta }}\Bigl \vert  (G_{N, \alpha})_{i+1,j+1} \Bigr \vert \vert x_{i+1} \vert \vert y_{j +1}\vert
 \le C_\alpha \sum_{(i,j)\in L_{1,\delta }}\vert i-j\vert ^{2\alpha-1} \vert x_{i+1} \vert \vert y_{j +1}\vert.$$
 Nous allons majorer cette somme pour $(i,j) \in [0,N_1 ]\times [2N_1,N-2N_1 ] =L_{1,\delta}$. Le 
 r\'esultat s'\'etendra sur tout $L_{2,\delta}$,$ L_{3,\delta}$, $L_{4,\delta}$ en utilisant les m\^{e}mes types de calculs.
 En remarquant que puisque $\vert i-j \vert \ge N\delta $ implique 
$$\left( \sum_{j=2N\delta }^{N-2N\delta } \vert i-j\vert ^{4\alpha-2} \right)^{1/2} 
< (N\delta)^{2\alpha-1} \left( \sum_{j=2N\delta }^{N-2N\delta } 1 \right)^{1/2} \le (N\delta )^{2\alpha_{1}-1/2}
$$
 et en se souvenant que 
 $$  \sum_{i=0}^{N\delta } \vert x_{i}\vert \le \left( \sum_{i=0}^{N\delta } \vert x_{i}\vert^2\right)^{1/2}
 \left( \sum_{i=0}^{N\delta } 1\right)^{1/2}$$
 On obtient 
 \begin{align*}
\sum_{(i,j)\in L_{1,\delta }}\vert i-j\vert ^{2\alpha-1} \vert x_{i+1} \vert \vert y_{j +1}\vert 
&= \sum_{i=0}^{N\delta } \vert x_{i}\vert \sum_{j=2N\delta }^{N-2N\delta } \vert i-j\vert ^{2\alpha-1} 
\vert y_{j} \vert \\
&\le \sum_{i=0}^{N\delta } \vert x_{i}\vert
\left( \sum_{j=2N\delta }^{N-2N\delta } \vert i-j\vert ^{4\alpha-2} \right)^{1/2} 
\left(  \sum_{j=2N\delta }^{N-2N\delta } \vert y_{j} \vert ^2\right)^{1/2}\\
&\le  \sum_{i=0}^{N\delta } \vert x_{i}\vert (N\delta )^{2\alpha-1/2}\\
& \le (N\delta )^{1/2} (N\delta )^{2\alpha-1/2} = (N\delta )^{2\alpha} =o(N^{2\alpha}).
\end{align*}

  En utilisant le lemme \ref{LEMME1} on obtient, pour $k=1$  
$$
  \Bigl \vert \sum _{(i,j) \in I_{1,\delta} } 
  (G_{N,\alpha})_{i+1,j+1} x_{j+1}  y_{i+1}\Bigr \vert 
 \le C_{\alpha}  \sum _{(i,j) \in I_1 } \vert i-j\vert ^{2\alpha-1} \vert x_{j+1}  y_{i+1}\vert $$
et il vient 
 $$
 \Bigl \vert \sum _{(i,j) \in I_{1,\delta } }
 (G_{N,\alpha})_{i+1,j+1} x_{j+1}  y_{i+1}\Bigr \vert 
 \le C_{\alpha}N^{2\alpha} \int_{-\delta }^\delta t^{2\alpha-1} dt =O\left( (N\delta )^{2\alpha}\right) = o(N^{2\alpha}).
 $$
 Les sommes portant sur $L_{2,\delta },$ $L_{3,\delta }$ et $L_{4,\delta }$ se traitent de m\^{e}me.\\ 
 Nous devons maintenant \'evaluer les quantit\'es 
 $$
 \Bigl \vert \sum _{(i,j) \in I_{k,\delta } }\left(T_{N}^{-1}(\varphi_{\alpha})\right)_{i+1,j+1} x_{j+1}  y_{i+1}\Bigr \vert
 $$
 pour $k\in \{1,2,3,4\}$.
 Des calculs prcis utilisant les thormes \ref{theoremepolypredi} et \ref{rappel} permettent d'obtenir que 
 $$
 \Bigl \vert \sum _{(i,j) \in I_{k,\delta } }\left(T_{N}^{-1}(\varphi_{\alpha})\right)_{i+1,j+1} x_{j+1}  y_{i+1}\Bigr \vert =o(N^{2\alpha})$$
(voir l'appendice). 
  Ce qui ach\`eve la d\'emonstration du lemme.
 \end{preuve}
   \subsection{D\'emonstration du th\'eor\`eme \ref{Principal} proprement dit}
    Rappelons tout d'abord le lemme (voir \cite{BoVi})
   \begin{lemme} \label{WIDOM}
 Soit $A_{N}=(a_{i,j})_{i,j=0}^{N-1}$ une matrice $N \times N$ \`a coefficients complexes. Soit
 $G_{N}$ l'op\'erateur int\'egral sur $L^2 [0,1]$ de noyau
 $$ g_{N}(x,y) = a_{[Nx], [Ny]}, \quad (x,y) \in (0,1)^2.$$
 Alors la norme de la matrice  $A_{N}$ et la norme de l'op\'erateur $G_{N}$ v\'erifient l'\'egalit\'e
 $\Vert A_{N}\Vert = N \Vert G_{N} \Vert.$
 \end{lemme}
Donnons nous un r\'eel $\mu$ v\'erifiant $1>\mu>1-\alpha$ et consid\'erons les op\'erateurs sur $L^2 (0,1)$ $G_{N,\alpha}^1, G_{N,\alpha}^2$, $(\tilde G_{N,\alpha})^1$
 et $(\tilde G_{N,\alpha})^2$ dont les noyaux $g_{N,\alpha}^1, g_{N,\alpha}^2$, $(\tilde g_{N,\alpha})^1$
 et $(\tilde g_{N,\alpha})^2$sont d\'efinis par  
 $$ g_{N,\alpha}^1 (x,y) =
\left \{
\begin{array}{ccc}
  N^{-2\alpha+1} G_{N,\alpha} ( [Nx],[Ny])  & \mathrm{si}   &  \vert x-y\vert >N^{\mu-1} \\
  0&  \mathrm{sinon} &   \end{array}
\right.
$$
  $$ g_{N,\alpha}^2 (x,y) =
\left \{
\begin{array}{ccc}
  N^{-2\alpha+1} G_{N,\alpha} ( [Nx],[Ny])  & \mathrm{si}   &  \vert x-y\vert <N^{\mu-1} \\
  0&  \mathrm{sinon} &   \end{array}
\right.
$$
 $$ (\tilde g_{N,\alpha})^1 (x,y) =
\left \{
\begin{array}{ccc}
   G_{\alpha} ( x,y)  & \mathrm{si}   &  \vert x-y\vert >N^{\mu-1} \\
  0&  \mathrm{sinon} &   \end{array}
\right.
$$
$$ (\tilde g_{N,\alpha})^2 (x,y) =
\left \{
\begin{array}{ccc}
   G_{\alpha} ( x,y)  & \mathrm{si}   &0<  \vert x-y\vert <N^{\mu-1} \\
  0&  \mathrm{sinon} &   \end{array}
\right.
$$
Dans la suite nous poserons $ x_{N}= \frac{[Nx]}{N}, y_{N}= \frac{[Ny]}{N}$, et supposerons $y>x$.
Il nous faut d'abord montrer que 
\begin{equation} \label{ROUTINE1}
\Vert G_{N,\alpha}^1 - (\tilde G_{N,\alpha})^1\Vert =o(1),
\end{equation}
c'est \`a dire que 
$\vert G_{\alpha}(x_{N},y_{N}) -G_{\alpha}(x,y)\vert =o(1)$ uniform\'ement pour 
$\vert x-y\vert > N^{\mu-1}$.\\
Dans un premier temps \'etudions,
$ \vert x^\alpha_{N}-x^\alpha \vert y^\alpha \int_{y}^1 \frac{(t-x)^{\alpha-1} (t-y)^{\alpha-1}}{t^{2\alpha}} dt $. 
On a facilement 
 \begin{align*}
  y^\alpha \int_{y}^1 \frac{(t-x)^{\alpha-1} (t-y)^{\alpha-1}}{t^{2\alpha}} dt
&\le  y (y-x)^{\alpha-1} \int_{y}^1 \frac{(t-y)^{\alpha-1}}{t^{\alpha+1}} dt \\
&\le  (N^{\mu-1})^{\alpha-1} (1-y)^\alpha \frac{1}{y\alpha} y =O((N^{\mu-1})^{\alpha-1}
\end{align*}
et finalement si $x_{N} \ge \frac{1}{N} $ le th\'eor\`eme des accroissements finis donne si $\mu>\alpha-1$
$$ \vert x_{N}^\alpha-x^\alpha \vert y^\alpha \int_{y}^1 \frac{(x-t)^{\alpha-1} (y-t)^{\alpha-1}}{t^{2\alpha}} dt 
= O(\frac{1}{N^{\alpha-1} }) \frac{1}{N} N^{(\mu-1)(\alpha-1)} =O( N^{(\mu-1) (\alpha-1)-\alpha}).$$
 Par contre si $x_{N}<\frac{1}{N} $ on a aussi $x_{N}=0$ et 
$$ \vert x_{N}-x \vert y^\alpha \int_{y}^1 \frac{(x-t)^{\alpha-1} (y-t)^{\alpha-1}}{t^{2\alpha}} dt 
= O(\frac{1}{N^{\alpha} })  N^{(\mu-1)(\alpha-1)} =O(N^{(\mu-1) (\alpha-1)-\alpha} )$$
On vrifie que si $\mu>1-\alpha$ alors $O(N^{(\mu-1) (\alpha-1)-\alpha} )=o(1)$.\\
Nous avons de m\^{e}me  
\begin{align*}
 \vert y^\alpha_{N}-y^\alpha \vert x_{N}^\alpha &\int_{y}^1 \frac{(x-t)^{\alpha-1} (y-t)^{\alpha-1}}{t^{2\alpha}} dt \\
 &\le  x_{N} \vert y^\alpha_{N}-y^\alpha \vert O(N^{\mu-1})^{\alpha-1} \frac{(1-y)^\alpha}{y\alpha} = O(N^{\mu-1})^{\alpha-1} N^{-\alpha} =o(1).
 \end{align*}
 Il nous faut ensuite consid\'erer 
\begin{align*}
x_{N}^\alpha y_{N}^\alpha
 \int_{y_{N}}^y \frac{(t-x)^{\alpha-1} (t-y)^{\alpha-1}}{t^{2\alpha}} dt 
&\le O( N^{(\mu-1) (\alpha-1)} )   \int_{y_{N}}^y (t-y)^{\alpha-1} dt 
\\
&= O(N^{(\mu-1) (\alpha-1)} N^{-\alpha} ) =o(1)
\end{align*}
avec l'hypoth\`ese faite sur $\mu$.
Consid\'erons maintenant la d\'ecomposition
\begin{align*}
& x_{N}^\alpha y_{N}^\alpha
 \int_{y}^1 \frac{(t-x)^{\alpha-1} (t-y)^{\alpha-1}-(t-x_{N})^{\alpha-1} (t-y_{N})^{\alpha-1}}{t^{2\alpha}} dt
  = \\ & x_{N}^\alpha y_{N}^\alpha
 \int_{y}^1 \frac{(t-x)^{\alpha-1} \left((t-y)^{\alpha-1}-(t-y_{N}\right))^{\alpha-1}}{t^{2\alpha}} dt\\
&+  x_{N}^\alpha y_{N}^\alpha
 \int_{y}^1 \frac{\left( (t-x)^{\alpha-1} -(t-x_{N})^{\alpha-1} \right)(t-y)^{\alpha-1}}{t^{2\alpha}} dt
 \\
 & =I_{1}+ I_{2}.
\end{align*}
On v\'erifie facilement que 
\begin{align*}
\vert I_{1}\vert &\le O(N^{(\mu-1)(\alpha-1)})   \int_{y}^1\vert (t-y)^{\alpha-1} -(t-y_{N})^{\alpha-1}\vert  dt\\
& \le O(N^{(\mu-1)(\alpha-1)})  \int_{y}^1 \left((t-y)^{\alpha-1}- (t-y_{N})^{\alpha-1} \right) dt\\
& \le  O\left(N^{(\mu-1)(\alpha-1)} \left( (1-y)^\alpha - (1-y_{N})^\alpha \right)\right)
\end{align*}
Si $\frac{1}{N}\le 1-y$ alors $\frac{1}{N}\le 1-y_{N}$ et 
$\left( (1-y)^\alpha - (1-y_{N})^\alpha \right) = \alpha c_{N}^{\alpha-1} \frac{1}{N}$ avec $c_{N}\ge \frac{1}{N}$.
Nous pouvons finalement conclure que $\vert I_{1}\vert =O(N^{(\mu-1)(\alpha-1)} N^{-\alpha})$ et 
puisque $\mu> \alpha-1$ on a  $(\mu-1)(\alpha-1) <\alpha$, c'est \`a dire que 
$\vert I_{1}\vert =o(1)$. \\
Si $\frac{1}{N}> 1-y$ alors $\frac{1}{N}\ge 1-y_{N}$ et 
$\left( (1-y)^\alpha - (1-y_{N})^\alpha \right) =O(N^{-\alpha})$ ce qui nous ram\`ene au cas pr\'ec\'edent.\\
Occupons nous maintenant de l'int\'egrale $I_{2}$. Nous pouvons \'ecrire, gr\^{a}ce au th\'eor\`eme des accroissements finis,  
\begin{align*} 
\vert I_{2}\vert &\le 
\int_{y}^1 \vert (t-x)^{\alpha-1} -(t-x_{N})^{\alpha-1} \vert (t-y)^{\alpha-1} dt \\
& \le \int_{y}^1 c_{N}^{\alpha-2}(t) \frac{1}{N} (\alpha-1)  (t-y)^{\alpha-1} dt
\end{align*}
avec $ N^{\mu-1} < y-x <t-x <c_{N}(t) <t-x_{N}$. 
Ce qui donne finalement 
$\vert I_{2}\vert = O(N^{(\mu-1)(\alpha-1)} ) =o(1)$.
En combinant les diverses majorations obtenues on obtient l'\'equation \ref{ROUTINE1}.\\
Il nous faut maintenant tudier 
$\Vert G^2_{N,\alpha}\Vert$. On montre comme 
dans \cite{BoVi} que pour estimer cette norme il suffit d'estimer $\Vert G^3_{N,\alpha}\Vert$ o 
$G^3_{N,\alpha}$ est la matrice dfinie par 

$$( G^3_{N,\alpha} )_{k+1,l+1}=
\left\{
\begin{array}{ccc}
 N^{-2\alpha+1} G_{\alpha}(k,l)  &  \mathrm{si} &  
 0<\vert k-l\vert \le N^\mu \\
 0 &  \mathrm{sinon}  &   \end{array}
\right.
$$
En se rapportant  la dmonstration du lemme 
\ref{Fonda} on comprend que 
$\Vert G^3_{N,\alpha}\Vert = o(N^{2\alpha})  N^{-2\alpha+1} = o(N)$
et donc 
$\Vert G^2_{N,\alpha}\Vert=o(1)$ en utilisant le lemme \ref{WIDOM}.\\
Enfin avec le lemme \ref{LEMME1} 
on obtient $\Vert (\tilde G^2_{N,\alpha})\Vert$ est
 majore par $O(\Vert K^2_\alpha\Vert)$ o 
 $K^2_\alpha$ est l'oprateur dans $L^2(0,1)$ 
 de noyau $(k_{\alpha}^2)$ dfini par 
\[ (k_{\alpha}^2) (x,y) =
\left\{
\begin{array}{ccc}
  \vert x-y\vert^{2\alpha-1} &  \mathrm {si} &  \vert x-y\vert < N^{\mu-1} \\
 0 &  \mathrm {sinon}   & \end{array}
\right.
\]
 On sait d'autre part que  $\Vert K^2_\alpha\Vert=o(1)$(voir \cite{BoVi}), ce qui achve de dmontrer la proprit.

   \section{D\'emonstration du corollaire \ref{encadrement}}

En utilisant le lemme \ref{LEMME1} et le thorme \ref{Noyau} on obtient 
\[ \Lambda_{N,\alpha} \le \frac{ N^{2\alpha}}{\Gamma (\alpha) c_{1}(1) } \frac{\Gamma (1-2\alpha)}
{\Gamma(1-\alpha)} \Vert K_{\alpha}\Vert.
\]
Comme l'on sait que 
$\Vert K_{\alpha}\Vert \le \frac{1}{\alpha}$ (voir \cite{BoVi}) on en d\'eduit la majoration 
\[ \Lambda_{N,\alpha} \le \frac{ N^{2\alpha}}{\Gamma (\alpha) c_{1}(1) } \frac{\Gamma (1-2\alpha)}
{\Gamma(1-\alpha)} \frac{1}{\alpha}.
\]
et donc la minoration de la plus petite valeur propre de 
$T_{N}\left( \vert 1-\chi \vert^{2\alpha} c\right)$
\\ 
D'autre part nous pouvons \'ecrire $\Vert \tilde G_{\alpha}\Vert \ge \Vert \tilde G_{\alpha} (\mathbf {1})\Vert$ si $\mathbf {1}$ dsigne la fonction constante gale  $1$.
En r\'eutilisant la minoration de $G_{\alpha}(x,y)$ utilis\'ee 
 dans \cite {RS1111} pour minorer $\Lambda_{N,\alpha}$ quand $\alpha\in ]\frac{1}{2},1[$ (paragraphe 7.1)
 ( savoir $G_{\alpha}(x,y)\ge x^\alpha y^\alpha (1-x)
 (1-y)$) ,
 minoration qui est toujours valable ici, il vient 
\begin{align*}
  \Vert \tilde G_{\alpha} (\textbf{1})\Vert &\ge \left( \int_{0}^1x^{2\alpha} (1-x)^2
  ( \int_{0}^1  y^\alpha (1-y) dy)dx \right)^{1/2} \\
 & \ge \left( \int_{0}^1x^{2\alpha} (1-x)^2
  ( \int_{0}^1  y^{2\alpha} (1-y)^2 dy)dx \right)^{1/2}\\
 &\ge  \int_{0}^1  x^{2\alpha} (1-x)^2 dy = \frac{6\Gamma (1+2\alpha)}{\Gamma (2\alpha+4)}.
\end{align*}
 
  \section{Dmonstration du thorme \ref{PROD}}
 Comme dans la dmonstration du lemme 
 \ref{Fonda} nous dfinissons les deux matrices 
 $G_{N,\alpha_1}$ et 
 $G_{N,\alpha_2}$ par 
 $$(G_{N,\alpha_1})_{k+1,l+1} =
 N^{2\alpha_1-1} \frac{1}{c_1(1)} G_{\alpha_1}(\frac{k}{N}, \frac{l}{N}) \quad \mathrm{si} \quad k\neq l
 \quad \mathrm{et} \quad (G_{N,\alpha_1})_{k+1,k+1}=0.$$
 $$(G_{N,\alpha_2})_{k+1,l+1} =
 N^{2\alpha_1-1} \frac{1}{c_1(1)} G_{\alpha_2}(\frac{k}{N}, \frac{l}{N}) \quad \mathrm{si} \quad k\neq l
 \quad \mathrm{et} \quad (G_{N,\alpha_2})_{k+1,k+1}=0.$$
 
 En utilisant le lemme \ref{Fonda} de la dmonstration du thorme \ref{Principal} il vient 
 \begin{align*}
 \Vert T_N^{-1}(\varphi_{\alpha_1})
 T_N^{-1} &(\varphi_{\alpha_2})  - G_{N,\alpha_1}
 G_{N,\alpha_2} \Vert = \\
 &= 
  \Vert T_N^{-1}(\varphi_{\alpha_1})
 T_N^{-1} (\varphi_{\alpha_2}) -
 T_N^{-1}(\varphi_{\alpha_1})   G_{N,\alpha_2}
 +T_N^{-1}(\varphi_{\alpha_1})   G_{N,\alpha_2}
 - G_{N,\alpha_1} G_{N,\alpha_2} \Vert \\
& \le \Vert T_N^{-1}(\varphi_{\alpha_1})
 \left( T_N^{-1} (\varphi_{\alpha_2}) -
   G_{N,\alpha_2}\right) \Vert 
   + \Vert  G_{N,\alpha_2} \left( 
   T_N^{-1}(\varphi_{\alpha_1}) - G_{N,\alpha_1}\right) \Vert\\
   &\le O(N^{2\alpha_1}) o(N^{2\alpha_2})
   +  O(N^{2\alpha_2}) o(N^{2\alpha_1}) 
   = o(N^{2\alpha_1+2\alpha_2}).
   \end{align*}
 Nous sommes donc ramen  valuer 
$\Vert  G_{N,\alpha_1} G_{N,\alpha_2} \Vert$.
Pour ce faire notons l'oprateur 
$G_{N,\alpha_1+\alpha_2}$ dfini sur $L^2(0,1)$ par
 $$ (x,y) \rightarrow 
  N^{2\alpha_1+2\alpha_2-1} \sum_{h\neq [Nx],
  h\neq [Ny]} (G_{N,\alpha_1})_{[Nx]+1,h} 
(G_{N,\alpha_2})_{h,[Ny]+1}.$$
Nous allons en fait montrer que 
\begin{equation}\label{BUT1}
\Vert G_{N,\alpha_1+\alpha_2} -
\tilde G_ {\alpha_1}\star \tilde G_{\alpha_2}\Vert = o(1).
\end{equation}
Reprenons la notation $x_N=\frac{[Nx]} {N}$
et $y_N=\frac{[Ny]} {N}$, supposons que $y>x$ et donnons nous un rel $\delta\in]0,1[$. 
Posons $J_{\delta,x_N} = ]x_N -\frac{[N \delta]}{N}, x_N +\frac{[N\delta]}{N}[$ 
et $J_{\delta,y_N} = ]y_N -\frac{[N \delta]}{N}, y_N +\frac{[N\delta]}{N}[$. Nous noterons 
  $N J_{\delta,x_N} = ][Nx] -[N \delta], 
  [Nx] +[N\delta][$ 
et $N J_{\delta,y_N} = ][Ny] -[N \delta], 
  [Ny] +[N\delta][$.
Considrons 
les oprateurs sur $L^2(0,1)$ 
$G_{N,\delta,\alpha_1+\alpha_2}^1$,
$G_{N,\delta,\alpha_1+\alpha_2}^2$, 
de noyaux respectifs 
$g_{N,\alpha_1+\alpha_2,\delta}^1$,
$g_{N,\alpha_1+\alpha_2,\delta}^2$, 
qui sont dfinis par
$$
g_{N,\alpha_1+\alpha_2,\delta}^1(x,y) =
\frac{1}{N} \sum_{h \in [0,N] \setminus (N J_{\delta,x_N}\cup N J_{\delta,y_N})}  
G_{\alpha_1} (x_N,\frac{h}{N}) 
G_{\alpha_2}(\frac{h}{N},y_N),$$
et 
$$
g_{N,\alpha_1+\alpha_2,\delta}^2(x,y) =
\frac{1}{N} \sum_{h\in N J_{\delta,x_N}\cup N J_{\delta,y_N}}  
G_{\alpha_1} (x_N,\frac{h}{N}) 
G_{\alpha_2}(\frac{h}{N},y_N).$$
Considrons tout d'abord la diffrence 
$$ 
D =g^1_{N,\alpha_1+\alpha_2,\delta}(x,y) 
- \int_{[0,1]\setminus (J_{\delta,x_N}\cup J_{\delta,y_N})} G_{\alpha_1}(x,t)G_{\alpha_2}(t,y) dt.
$$
Nous pouvons crire 
\begin{align*} 
D &= \sum_{h\in [0,N-1] \setminus (N J_{\delta,x_N}\cup N J_{\delta,y_N})} \int_{h/N}^{(h+1)/N} 
\left(G_{\alpha_1} (x_N,\frac{h}{N}) 
G_{\alpha_2}(\frac{h}{N},y_N) -
 G_{\alpha_1}(x,t)G_{\alpha_2}(t,y) \right) dt\\
 &=  \sum_{h\in [0,N-1] \setminus (N J_{\delta,x_N}\cup N J_{\delta,y_N})} \int_{h/N}^{(h+1)/N} 
\left(G_{\alpha_1} (x_N,\frac{h}{N}) -G_{\alpha_1}(x,t) \right) G_{\alpha_2}(\frac{h}{N},y_N) dt 
\\
&+ \sum_{h\in [0,N-1] \setminus (N J_{\delta,x_N}\cup N J_{\delta,y_N})} \int_{h/N}^{(h+1)/N} 
G_{\alpha_1} (x,t) \left( G_{\alpha_2}(\frac{h}{N},y_N)
- G_{\alpha_2} (t,y) \right)dt.
\end{align*}
En remarquant que si 
$t\in [\frac{h}{N}, \frac{h+1}{N}]$
alors $[Nt]=h$ nous pouvons crire,
en utilisant les rsultats acquis dans la dmonstration du thorme \ref{Principal}, 
que 
$\vert G_{\alpha_1} (x_N,\frac{h}{N}) -G_{\alpha_1}(x,t)\vert =o(1)$
et 
$\vert G_{\alpha_2} (\frac{h}{N},y_N) -G_{\alpha_2}(t,y)\vert =o(1)$
 uniformment en $x,y,t$ (en effet si 
 $\vert x-t \vert >\delta$ alors 
 $\vert x-t \vert >N^{\mu-1}$ si $N$ assez grand) .
D'o si $\epsilon>0$ assez petit et $N$ suffisamment grand 
\begin{align*}
&\sum_{h\in [0,N-1] \setminus (N J_{\delta,x_N}\cup N J_{\delta,y_N})} \int_{h/N}^{(h+1)/N} 
\left(G_{\alpha_1} (x_N,\frac{h}{N}) -G_{\alpha_1}(x,t) \right) G_{\alpha_2}(\frac{h}{N},y_N)  dt\\
 &\le\sum_{h\in [0,N] \setminus (N J_{\delta,x_N}\cup N J_{\delta,y_N})}\epsilon \int_{h/N}^{(h+1)/N} 
G_{\alpha_2}(t,y) dt +\epsilon.
\end{align*}
uniformment en $x$ et $y$.\\
D'autre part , en utilisant  le lemme \ref{LEMME1}, on a :
$$ \int_0^1 G_{\alpha_2}(t,y) dt 
=O\left( \int_0^1 \vert t-y\vert ^{2\alpha_2-1} dt \right)
=O(1).$$
En traitant de m\^{e}me le terme
$$  \sum_{h\in [0,N-1] \setminus (N J_{\delta,x_N}\cup N J_{\delta,y_N})} \int_{h/N}^{(h+1)/N} 
G_{\alpha_1} (x,t) \left( G_{\alpha_2}(\frac{h}{N},y_N)
- G_{\alpha_2} (t,y)\right) dt$$
on obtient $\vert D\vert =o(1)$ uniformment en $x$ et $y$. \\  
D'autre part il vient, toujours avec le lemme \ref{LEMME1} 
$$
\vert g_{N,\alpha_1+\alpha_2,\delta}^2(x,y)\vert 
\le \frac{1}{N}  \sum_{h\in N J_{\delta,x_N}\cup N J_{\delta,y_N}} \Bigr \vert x_N - \frac{h}{N} \Bigl\vert ^{2\alpha_1-1} \Bigr \vert \frac{h}{N} -y_N\Bigl \vert ^{2\alpha_2-1}. 
$$
En utilisant la monotonie 
de la fonction $t\rightarrow \vert x-t\vert ^{2\alpha_1-1}  \vert x-t\vert ^{2\alpha_2-1} $ comme dans \cite{Ramb10}, on obtient 
$$ \sum_{h\in N J_{\delta,x_N}\cup N J_{\delta,y_N})}  \vert x_N - \frac{h}{N} \vert ^{2\alpha_1-1} \vert \frac{h}{N} -y_N\vert ^{2\alpha_2-1}
\sim \int_{(J_{\delta,x_N}\cup  J_{\delta,y_N})}
\vert x-t\vert ^{2\alpha_1-1}  \vert x-t\vert ^{2\alpha_2-1} dt =O(\delta).$$
En faisant maintenant tendre $\delta$ vers zro et en 
utilisant encore une fois le lemme \ref{LEMME1} 
pour obtenir la convergence de l'intgrale on obtient (\ref{BUT1}).
\section{Dmonstration du thorme 
\ref{encadrement2}}
On a obtenu dans \cite{Ramb10} l'encadrement suivant   \begin{lemme}\label{DUO}
           Si $0<\alpha_{1}\le\alpha_{2}<\frac{1}{2}$ et $x\neq y$ nous avons 
           \[ \vert x-y\vert ^{2\alpha_{1}+2\alpha_{2}-1} \le \int_{0}^1 \vert x-t \vert ^{2\alpha_{1}-1}
           \vert y-t \vert ^{2\alpha_{2}-1} dt \le H_{\alpha_{1},\alpha_{2} }\vert x-y \vert ^{2\alpha_{1}+2\alpha_2-1}\]
          avec
          $ H_{\alpha_1,\alpha_2}= \frac{1}{\alpha_{1}}+ \frac{1}{\alpha_{2}}$.
                   \end{lemme}
Ce lemme fournit une majoration immdiate de 
$\Vert G_{N,\alpha_1}\star G_{N,\alpha_2} \Vert.$
Pour minorer cette norme on peut remarquer qu'elle 
est suprieure  
$\Vert (G_{N,\alpha_1}\star G_{N,\alpha_2})
 \mathbf 1\Vert$ que l'on peut minorer, toujours en utilisant 
les minorations de \cite{RS1111}, par 
\begin{align*}
&\left( \int_0^1 \left(\int_0^1 (G_{\alpha_1}\star 
G_{\alpha_2}) (x,y)dy\right)^2dx\right)^{1/2}\\
&\ge \left( \int_0^1 \left(\int_0^1
\left( \int_0^1 x^{\alpha_1} (1-x)  
t^{\alpha_1+\alpha_2} (1-t)^{2}
y^{\alpha_2} (1-y) dt\right) dy\right)^2dx\right)^{1/2}\\
&\ge \left( \int_0^1 x^{2\alpha_1} (1-x)^2 dx \right)^{1/2}
\left(\int_0^1 y^{\alpha_2} (1-y) dx \right)
\left(\int_0^1 t^{\alpha_1+\alpha_2} (1-t)^2 dx \right)\\
&\ge \left( \int_0^1 x^{2\alpha_1} (1-x)^2 dx \right)
\left(\int_0^1 y^{2\alpha_2} (1-y)^2 dx \right)
\left(\int_0^1 t^{\alpha_1+\alpha_2} (1-t)^2 dx \right)
\end{align*}
ce qui achve la dmonstration.

\section{Appendice}
Nous devons majorer les quantits 
$$
 \Bigl \vert \sum _{(i,j) \in I_{k,\delta } }\left(T_{N}^{-1}(\varphi_{\alpha})\right)_{i+1,j+1} x_{j+1}  y_{i+1}\Bigr \vert $$
 pour $k \in \{1,2,3,4\}$.
Nous allons nous concentrer sur 
 $$
 \Bigl \vert \sum _{(i,j) \in I_{1,\delta } }\left(T_{N}^{-1}(\varphi_{\alpha})\right)_{i+1,j+1} x_{j+1}  y_{i+1}\Bigr \vert $$
 Pour cel\`a nous devons majorer avec pr\'ecision les quantit\'es 
 $\left(T_{N}^{-1}(\varphi_{\alpha})\right)_{i+1,j+1}$ pour $(i,j)\in I_{1,\delta}$. 
Utilisons encore la formule (\ref{GOBSEM}). En supposant 
  $i\le j$, le th\'eor\`eme \ref {theoremepolypredi} permet d'\'ecrire 
  $$ \sum_{u=0}^i \overline{\gamma^{(\alpha)}_{i-u}} \gamma^{(\alpha)}_{j-u} 
  = \left(\sum_{u=0}^i \overline{\beta^{(\alpha)}_{i-u}} \beta^{(\alpha)}_{j-u}\right) \left(1+o(1)\right). $$
  Notons $k_{0}$ d\'esigne un entier tel que $\beta^{(\alpha)}_{k}$ puisse \^{e}tre 
  remplac\'e par son asymptotique pour $k\ge k_{0}.$ Nous sommes amen\'es \`a distinguer quatre cas.
  \begin{itemize}
  \item [$\bullet$] 
 Si $j\ge i\ge k_0$ et $0 \le j-i\le k_{0}$ on crit 
   $$
  \sum_{u=0}^i \overline{\beta^{(\alpha)}_{i-u}} \beta^{(\alpha)}_{j-u}  = 
 \sum _{u=i-k_0+1}^i \overline{\beta^{(\alpha)}_{i-u}} \beta^{(\alpha)}_{j-u} + \sum_{u=0}^{i-k_0} \overline{\beta^{(\alpha)}_{i-u}} \beta^{(\alpha)}_{j-u}.
 $$
En posant $M_{1}= \displaystyle{ \sum_{0\le h_{1}\le k_{0}, 0\le h_{2}\le 2 k_{0}} 
\Bigl \vert\overline{\beta^{(\alpha)}_{i-u}} \beta^{(\alpha)}_{j-u}\Bigr \vert}$ on 
obtient 
 \begin{eqnarray*}
\Bigl \vert \sum _{u=i-k_0+1}^i \overline{\beta^{(\alpha)}_{i-u}} \beta^{(\alpha)}_{j-u}\Bigr\vert &\le& 
M_{1}= (M_{1} k_0^{1-\alpha}) 
 k_0^{\alpha-1} \\
&\le& (M_{1} k_0^{1-\alpha}) (j-i)^{\alpha-1}\le 
 (M_{1} k_0^{1-\alpha}) \vert j-i\vert^{2\alpha-1}
 \end{eqnarray*}
 et avec le lemme \ref{DUO} 
 $$
 \Bigl \vert \sum_{u=0}^{i-k_0} \overline{\beta^{(\alpha)}_{i-u}} \beta^{(\alpha)}_{j-u}\Bigr\vert 
  \sim
\frac{N^{2\alpha-1}}{\Gamma^2(\alpha) c_1(1)} \int _{0}^x (x-t)^{\alpha-1} (y-t)^{\alpha-1} dt \le
H_{\alpha_{1},\alpha_{2}} \vert j-i\vert ^{2\alpha-1}$$
en posant $ x = \frac{i}{N}$ et 
 $ y = \frac{j}{N}$.
 \item [$\bullet\bullet$] 
 Si $0\le i< k_0$ et $0 \le j-i\le k_{0}$ on peut alors crire, en remarquant que  
 $$
 \Bigl \vert \sum_{u=0}^i \overline{\beta^{(\alpha)}_{i-u}} \beta^{(\alpha)}_{j-u} \Bigr\vert  \le 
 M_{1}
 $$
et comme pr\'ec\'edemment 
$$
M_1 \le 
 (M_1 k_0^{1-\alpha}) \vert j-i\vert^{2\alpha-1}.$$

  \item [$\bullet\bullet\bullet$] 
 Si $j\ge i\ge k_0$ et $ j-i\ge k_{0}$ on crit 
   $$
  \sum_{u=0}^i \overline{\beta^{(\alpha)}_{i-u}} \beta^{(\alpha)}_{j-u}  = 
 \sum _{u=i-k_0+1}^i \overline{\beta^{(\alpha)}_{i-u}} \beta^{(\alpha)}_{j-u} + \sum_{u=0}^{i-k_0} \overline{\beta^{(\alpha)}_{i-u}} \beta^{(\alpha)}_{j-u}.
 $$
Si $M_2 = \max_{0\le h\le k_{0}}  \vert \beta_{h}^{(\alpha)}\vert $ nous pouvons \'ecrire 
 \begin{align*}
\Bigl \vert \sum _{u=i-k_0+1}^i \overline{\beta^{(\alpha)}_{i-u}} \beta^{(\alpha)}_{j-u}\Bigr\vert &\le
M_2 \sum _{u=i-k_0+1}^i \vert j-u\vert^{\alpha-1}\\
&\le M_2 k_0\vert j-i\vert^{\alpha-1} 
\le M_2 k_0\vert j-i\vert^{2\alpha-1},
\end{align*}
et de plus, toujours avec le lemme \ref{DUO} 
$$ \sum_{u=0}^{i-k_0} \overline{\beta^{(\alpha)}_{i-u}} \beta^{(\alpha)}_{j-u} \sim
\frac{ N^{2\alpha-1}}{\Gamma^2(\alpha) c_1(1)} \int _{0}^x (x-t)^{\alpha-1} (y-t)^{\alpha-1} dt \le 
H_{\alpha_{1},\alpha_{2}} \vert j-i\vert^{2\alpha-1}$$
par des calculs d\'ej\`a vu et toujours en posant en posant $ x = \frac{i}{N}$ et 
 $ y = \frac{j}{N}$.
 \item[$\bullet\bullet\bullet\bullet$]
 Si $0\le i< k_0$ et $ j-i\ge k_{0}$.
On peut alors \'ecrire 
 $$ \Bigl \vert _{u=0}^i\overline{\beta_{i-u} ^{(\alpha)}}
 \beta_{j-u} ^{(\alpha)}\Bigr \vert\sim \frac{\vert j-i\vert ^{\alpha-1} }{\Gamma (\alpha) c_{1}(1)} 
 \sum_{v=0}^{k_{0}} \vert \beta_{v}^{(\alpha)} \vert 
 .$$
  \end{itemize}
 On obtient  finalement 
 \begin{align*}  \Bigl \vert \sum _{(i,j) \in I_{1,\delta } }
&\left( \sum_{u=0}^i \overline{\gamma^{(\alpha)}_{i-u}} \gamma^{(\alpha)}_{j-u} \right) x_{i+1} y_{j+1}\Bigr \vert 
 \\ & \le O\left(N^{2\alpha} \int_{-\delta }^\delta t^{2\alpha-1} dt\right) =O\left( (N\delta )^{2\alpha}\right) = o(N^{2\alpha}).
 \end{align*}
Enfin le th\'eor\`eme \ref{rappel} permet d'\'ecrire, 
toujours si $i\le j$  
\begin{align*}
 \sum_{v=0}^{i} \gamma_{v+N-j}^{(\alpha)}\overline{\gamma_{v+N-i}^{(\alpha)} }
&\sim  \sum_{v=0}^{i}\frac{ \beta_{j-v}^{(\alpha+1)}}{N}
\frac{\overline{\beta_{i-v}^{(\alpha+1)} }}{N}\\
& = \frac{O((N\delta )^{2\alpha+1)}}{N^2} =o(N^{2\alpha-1}).
\end{align*}

On obtient alors 
$$
  \Bigl \vert \sum _{(i,j) \in I_{1,\delta } }
 \left(\sum_{u=N-j}^{N-j+i} \gamma_{u}^{(\alpha)}\overline{\gamma_{u+j-i}^{(\alpha)} }\right) x_{i+1}y_{j+1}\Bigr \vert 
\le o(N^{2\alpha-1}) \sum _{(i,j)\in I_{1,\delta } }
 \vert x_{i+1} y_{j+1}\vert= o(N^{2\alpha-1}) 
$$
puisque $\displaystyle{\sum _{(i,j)\in I_{1,\delta } }
 \vert x_{i+1} y_{j+1}\vert} \le 1.$\\
La majoration de la somme sur $I_{2,\delta }$ se d\'eduit de ce r\'esultat en utilisant les sym\'etries de la matrice $T_{N}^{-1}\left( \varphi_{\alpha}\right)$.
Les m\^{e}mes mthodes que pour $I_{1,\delta}$ donnent le rsultat sur  $I_{3,\delta }$ puis  $I_{4,\delta }$par symtrie.

  \bibliography{Toeplitzdeux}

\end{document}